\documentclass[11pt]{article}
\usepackage[T1]{fontenc}
\usepackage[utf8]{inputenc}
\usepackage[english]{babel}
\usepackage{geometry}
\usepackage{comment}
\usepackage{amssymb,amsmath,amsthm,mathtools,stmaryrd}
\usepackage{xcolor}
\usepackage{authblk}
\usepackage[colorlinks=true,citecolor=red,linkcolor=blue]{hyperref}

\newtheorem{theorem}{Theorem}[section]
\newtheorem{proposition}[theorem]{Proposition}
\newtheorem{lemma}[theorem]{Lemma}
\newtheorem{corollary}[theorem]{Corollary}
\theoremstyle{definition}
\newtheorem{assumption}[theorem]{Assumption}
\theoremstyle{remark}
\newtheorem{remark}[theorem]{Remark}
\numberwithin{equation}{section}

\newcommand{\trp}{\mathtt t^+}
\newcommand{\trm}{\mathtt t^-}

\providecommand{\keywords}[1]{\small\textbf{\textit{Keywords:}} #1}

\title{Carleman estimates and $\phi$-null controllability for a fully discrete fourth-order parabolic equation}
\author[1]{P\'ia Bizama}
\author[2]{V\'ictor Hern\'andez-Santamar\'ia}
\author[3]{Ariel A. P\'erez}
\author[4]{Jose Solarte-Ballesteros}
\affil[1]{Departamento de Matem\'atica, Universidad del B\'io-B\'io, Concepci\'on, Chile. \texttt{pbizama@ubiobio.cl}}
\affil[2]{Departamento de Matem\'aticas, Facultad de Ciencias
Universidad Nacional Aut\'onoma de M\'exico, M\'exico. \texttt{victor.santamaria@ciencias.unam.mx}}
\affil[3]{Departamento de Ingenier\'ia Matem\'atica, Universidad de Concepci\'on, Concepci\'on, Chile. \texttt{arperez@udec.cl}}
\affil[4]{Departmento de Matem\'atica, Universidad del B\'io-B\'io, Concepci\'on, Chile. \texttt{jsolarte@ubiobio.cl}}
\date{\today}

\begin{document}
\maketitle

\begin{abstract}
We consider a fully discrete approximation of a one-dimensional fourth-order parabolic equation with clamped boundary conditions, obtained by centered finite differences in space and the implicit Euler method in time, and controlled by two boundary controls acting at the left endpoint. We prove a Carleman estimate for the corresponding discrete adjoint operator, in which the large Carleman parameter is connected to the space and time discretization parameters $h$ and $\Delta t$. As a consequence, we obtain a relaxed observability inequality and a $\phi(h)$-null controllability result, in which the controls are bounded uniformly with respect to $h$ and $\Delta t$, and the norm of the solution at the final time is of order $\sqrt{\phi(h)}$, with $\phi(h)$ exponentially small, provided $\Delta t\le h^{\mu}$ for some $\mu\ge1$.
\end{abstract}

\keywords{Higher-order parabolic equations; Uniform controllability; Boundary observability, Relaxed observability; Penalized HUM.}\\
\textbf{Mathematics Subject Classification (2020):}  93B05, 93B07, 35K52, 65M06, 93C55.

\section{Introduction}\label{sec:intro}
\subsection{Motivation and background}

Fourth-order parabolic equations appear in the modelling of phase separation, through the Cahn--Hilliard equation \cite{cahn-1958}, and of flame-front propagation and instabilities, through the Kuramoto--Sivashinsky equation \cite{kuramoto-tsuzuki-1975,kuramoto-tsuzuki-1976,sivashinsky-1977}. Their controllability properties have been studied extensively in the continuous setting. In one space dimension, the boundary null controllability of the linear Kuramoto--Sivashinsky equation was established in \cite{cerpa-2010}, the case of a single boundary control was analysed in \cite{CGM:2017}, and the local exact controllability to the trajectories of the nonlinear equation with boundary controls was proved in \cite{cerpa-mercado} by means of a Carleman estimate. Observability and null controllability for one-dimensional fourth-order parabolic equations were also obtained in \cite{zhou-2012}. In several space dimensions, the case of interior controls was treated in \cite{guerrero-2019} by a global Carleman estimate, and in \cite{LRR:2020} by a spectral inequality for the clamped bi-Laplace operator. Boundary null controllability with two boundary controls follows from \cite{guerrero-2019} by extending the domain, whereas controllability with fewer boundary controls is largely open; it has been established for the Kuramoto--Sivashinsky equation on a rectangle \cite{takahashi-2017} and on cylindrical domains \cite{HSM:2025}, by means of the moment method combined, in the latter, with the Lebeau--Robbiano strategy. All these results follow the classical duality approach, see \cite{lions}, in which null controllability is equivalent to an observability inequality for the adjoint system, the latter being obtained from a global Carleman estimate in the spirit of \cite{FI,l-r}, or by the moment method \cite{f-r}.

From the numerical point of view, a natural question is whether these properties are preserved, uniformly with respect to the discretization parameters, when the equation is approximated. In general they are not, since even for the heat equation discrete observability inequalities may fail or have constants that blow up as the mesh is refined, see \cite{zuazua:2005} and the fully discrete counterexample in \cite{AAP-2024}. This lack of uniformity is closely related to the difficulties encountered in the numerical computation of null controls for parabolic problems, see for instance \cite{c-g-l,m-z,fc-m,fc-m2}. For space semi-discrete parabolic equations, uniform controllability properties were first obtained in \cite{LT:2006} by spectral arguments, and in \cite{BHL10a,boyer-sicon} by means of discrete Carleman estimates. Since an exact discrete null control cannot be expected to be uniformly bounded in general, the relevant notion is that of $\phi$-null controllability, in which one seeks controls that are bounded uniformly with respect to the discretization parameters and drive the solution to a state whose norm is of order $\sqrt{\phi}$, where $\phi$ is a prescribed function of the discretization parameters that tends to zero, typically at an exponential rate. This approach, based on a relaxed observability inequality and on the penalized Hilbert Uniqueness Method, is described in detail in \cite{boyer-canum}.

Discrete Carleman estimates and the resulting $\phi$-null controllability results have since been developed in several directions, including semilinear semi-discrete equations \cite{boyer-2014}, discontinuous diffusion coefficients \cite{Thuy}, abstract semi-discrete parabolic systems in Banach spaces \cite{Thuy2}, and coupled parabolic systems with boundary controls \cite{allonsius,allonsius2}. Space/time discretized parabolic equations were considered in \cite{boyer2011}. Time-discrete parabolic equations, for which the relevant notion is $\phi(\Delta t)$-null controllability, were treated in \cite{boyer-hernandez}, and a simplified time-discrete stabilized Kuramoto--Sivashinsky system was studied in \cite{HS:2023}. Carleman estimates for fully discrete parabolic operators, with both $h$ and $\Delta t$ connected to the Carleman parameters, were obtained for the one-dimensional heat equation with Dirichlet boundary conditions in \cite{GC-HS-2021}, with dynamic boundary conditions in \cite{LMPZ:2023}, and in arbitrary dimension in \cite{AAP-2024}. Beyond second-order parabolic operators, a fully discrete Carleman estimate for the third-order operator $\partial_t+\partial_x^3$ was recently established in \cite{kumar:2026}, where the $\phi$-null controllability of a fully discrete linear Korteweg--de Vries equation with a single boundary control is proved. All these works derive the Carleman estimate directly at the discrete level. A different route was proposed in \cite{LW:2026}, where an abstract transfer theorem derives fully discrete Carleman estimates for backward Euler schemes from a continuous Carleman inequality, through weighted error estimates for a shifted continuous auxiliary problem, and it is verified for the heat equation discretized by $P_1$ finite elements and by Cartesian finite differences with interior observation.

For higher-order operators, the discrete theory is much less developed. A Carleman estimate for a space semi-discrete fourth-order parabolic operator and the corresponding boundary $\phi$-null controllability result were obtained in \cite{CLTP:2022}; see also \cite{WZ:2024} for a stochastic semi-discrete fourth-order parabolic equation. A key ingredient in all these works is a precise description of the effect of discrete difference and average operators on the Carleman weight functions; for fourth-order operators these operators have to be applied up to order four, which makes the computations considerably heavier than in the second-order case. Such a description was obtained in one dimension in \cite{LOP:2021,LMPZ:2023} and in arbitrary dimension, including the time-discrete operator, in \cite{AAP-2024}.

To the best of our knowledge, no Carleman estimate with boundary observation is available for a fully discrete fourth-order parabolic operator. The purpose of this paper is to fill this gap in the one-dimensional case with boundary controls, and to deduce a $\phi$-null controllability result with controls bounded uniformly with respect to both $h$ and $\Delta t$.

\subsection{The discrete setting}\label{subsec:setting}

Let $T>0$, $M,N\in\mathbb N$, $h:=1/(M+1)$ and $\Delta t:=T/N$. In space we use the primal mesh, its boundary and its closure, and the set of ghost points with the corresponding extended mesh,
$$
\mathcal{M}:=\{ih\,;\,i\in\llbracket 1,M\rrbracket\},\quad \partial\mathcal{M}:=\{0,1\},\quad \overline{\mathcal M}:=\mathcal{M}\cup\partial\mathcal{M},\quad \partial\overline{\mathcal M}:=\{-h,1+h\},\quad \widehat{\mathcal M}:=\overline{\mathcal M}\cup\partial\overline{\mathcal M},
$$
the dual mesh and its extension,
$$
\mathcal M^{\ast}:=\{(i+\tfrac12)h\,;\,i\in\llbracket 0,M\rrbracket\},\qquad \partial\widehat{\mathcal M}^{\ast}:=\{-\tfrac h2,1+\tfrac h2\},\qquad \widehat{\mathcal M}^{\ast}:=\mathcal M^{\ast}\cup\partial\widehat{\mathcal M}^{\ast} .
$$
For a function $u$ defined on the nodes, the difference and average operators are
$$
D_hu(x):=\frac{u(x+\frac h2)-u(x-\frac h2)}{h},\qquad A_hu(x):=\frac{u(x+\frac h2)+u(x-\frac h2)}{2};
$$
they map functions on $\widehat{\mathcal M}$ to functions on $\widehat{\mathcal M}^{\ast}$, and functions on $\widehat{\mathcal M}^{\ast}$ to functions on $\overline{\mathcal M}$. In particular $D_h^2u$ is defined on $\overline{\mathcal M}$, $D_h^3u$ on $\mathcal M^{\ast}$ and $D_h^4u$ on $\mathcal{M}$ for $u$ defined on $\widehat{\mathcal M}$. In time we use the primal and dual meshes
$$
\mathcal{N}:=\{k\Delta t\,;\,k\in\llbracket1,N\rrbracket\},\quad \overline{\mathcal N}:=\mathcal{N}\cup\{0\},\quad
\mathcal{N}^{\ast}:=\{(k-\tfrac12)\Delta t\,;\,k\in\llbracket1,N\rrbracket\},\quad \overline{\mathcal N}^{\ast}:=\mathcal{N}^{\ast}\cup\{T+\tfrac{\Delta t}{2}\},
$$
the translations $\mathtt t^{\pm}f(t):=f(t\pm\frac{\Delta t}{2})$ and the discrete time derivative $D_tf:=(\trp f-\trm f)/\Delta t$, which maps functions on $\overline{\mathcal N}$ (resp. $\overline{\mathcal N}^{\ast}$) to functions on $\mathcal{N}^{\ast}$ (resp. $\mathcal{N}$). Discrete integrals and norms are defined in Section \ref{subsec:space}.

\subsection{The control problem}

Given $y_0\in L^2_h(\mathcal{M})$ and a bounded potential $a$ defined on $\mathcal{M}\times\mathcal{N}$, we set $\gamma:=\|a\|_{L^\infty_h}:=\max_{\mathcal{M}\times\mathcal{N}}|a|$ and consider the fully discrete controlled system
\begin{equation}\label{sys:primal}
\begin{cases}
D_ty+D_h^4\trp y+\trp(a\,y)=0, & \text{in }\mathcal{M}\times\mathcal{N}^{\ast},\\
\trp y(0,t)=u_1(t),\quad \trp y(1,t)=0, & t\in\mathcal{N}^{\ast},\\
D_h\trp y(-\tfrac h2,t)=u_2(t),\quad D_h\trp y(1+\tfrac h2,t)=0, & t\in\mathcal{N}^{\ast},\\
y(\cdot,0)=y_0, & \text{in }\mathcal{M},
\end{cases}
\end{equation}
where $y$ is defined on $\widehat{\mathcal M}\times\overline{\mathcal N}$ and the controls $u_1,u_2$ are defined on $\mathcal{N}^{\ast}$. System \eqref{sys:primal} is the standard finite-difference discretization in space of
$$
\partial_ty+\partial_x^4y+ay=0\ \text{ in }(0,1)\times(0,T),\qquad y(0,\cdot)=u_1,\ \partial_xy(0,\cdot)=u_2,\qquad y(1,\cdot)=\partial_xy(1,\cdot)=0,
$$
combined with the implicit Euler method in time; the ghost points of $\partial\overline{\mathcal M}$ carry the discrete Neumann conditions. The corresponding adjoint system reads
\begin{equation}\label{sys:adjoint}
\begin{cases}
-D_tq+D_h^4\trm q+a\,\trm q=0, & \text{in }\mathcal{M}\times\mathcal{N},\\
\trm q=0, & \text{on }\partial\mathcal{M}\times\mathcal{N},\\
D_h\trm q=0, & \text{on }\partial\widehat{\mathcal M}^{\ast}\times\mathcal{N},\\
q(\cdot,T+\tfrac{\Delta t}{2})=q_T, & \text{in }\mathcal{M},
\end{cases}
\end{equation}
where $q$ is defined on $\widehat{\mathcal M}\times\overline{\mathcal N}^{\ast}$, and $q_T\in L^2_h(\mathcal{M})$ is extended by zero to $\partial\mathcal{M}\cup\partial\overline{\mathcal M}$. The observed quantities are $D_h^2q(0,t)$ and $D_h^3q(\tfrac h2,t)$, $t\in\mathcal{N}^{\ast}$, which approximate $\partial_x^2q(0,t)$ and $\partial_x^3q(0,t)$; they are the quantities dual to $u_2$ and $u_1$, respectively, in the discrete duality identity of Lemma \ref{lem:duality}.

\subsection{Main results}

Let $\varphi$ satisfy Assumption \ref{ass:psi} below, which is the assumption made in the semi-discrete setting \cite[(7)--(8)]{CLTP:2022}. For $0<\delta\le\frac12$ and $\tau\ge1$ we set
$$
\theta(t):=\frac{1}{(t+\delta T)(T+\delta T-t)},\qquad
s(t):=\tau\theta(t),\qquad r:=e^{-s\varphi},\quad \rho:=r^{-1}.
$$
These are the weights of \cite{CLTP:2022}, with the roles of $r$ and $\rho$ interchanged, so that the estimates of discrete operators applied to $r$ and $\rho$ obtained there can be used directly. Our first result is the following Carleman estimate, where $\mathcal P^\star q:=-D_tq+D_h^4\trm q$ and the meshes $Q$, $Q_{\ast}$, $Q^{\ast}$, $\overline{Q}_{\ast}$ are defined in \eqref{eq:Qmeshes}.

\begin{theorem}\label{thm:intro:carleman}
There exist $C>0$, $\tau_0\ge1$, $h_0>0$ and $\varepsilon_0\in(0,1]$ such that, for all $T>0$, $\delta\in(0,\frac12]$, $\tau\ge\tau_0(T+T^2)$, $h\in(0,h_0]$ and $\Delta t>0$ satisfying
\begin{equation}\label{eq:intro:conditions}
\frac{\tau h}{\delta T^{2}}\le\varepsilon_0\qquad\text{and}\qquad \frac{\tau^{5}\Delta t}{\delta^{5}T^{10}}\le\varepsilon_0,
\end{equation}
and for all $q$ defined on $\widehat{\mathcal M}\times\overline{\mathcal N}^{\ast}$ and satisfying the boundary conditions of \eqref{sys:adjoint} at every $t\in\overline{\mathcal N}^{\ast}$, we have
\begin{equation}\label{eq:intro:carleman}
\begin{aligned}
&\int_{Q_{\ast}}s^{7}r^{2}|q|^{2}+\int_{Q^{\ast}}s^{5}r^{2}|D_hq|^{2}+\int_{\overline{Q}_{\ast}}s^{3}r^{2}|D_h^{2}q|^{2}+\int_{Q^{\ast}}s\,r^{2}|D_h^{3}q|^{2}\\
&\qquad+\int_{Q}\trm(s^{-1}r^{2})|D_tq|^{2}+\int_{Q_{\ast}}s^{-1}r^{2}|D_h^{4}q|^{2}\\
&\le C\bigg(\int_{Q}(\trm r)^{2}|\mathcal P^\star q|^{2}+\int_{\mathcal{N}^{\ast}}s^{3}r^{2}|D_h^{2}q|^{2}\Big|_{x=0}+\int_{\mathcal{N}^{\ast}}s\,r^{2}\big|_{x=0}|D_h^{3}q|^{2}\Big|_{x=h/2}+h^{-4}\int_{\mathcal{M}\times\partial\mathcal{N}}|\trp(rq)|^{2}\bigg).
\end{aligned}
\end{equation}
\end{theorem}

Compared with the second-order case \cite{GC-HS-2021,LMPZ:2023,AAP-2024}, where the time-discretization condition reads $\tau^{4}\Delta t(\delta^{4}T^{6})^{-1}\le\varepsilon_0$, the second condition in \eqref{eq:intro:conditions} involves one more power of $\tau/\delta$. It comes from the term of order $\Delta t\,s^{2}|D_hD_tz|^{2}$, with the unfavourable sign, produced by the second-order part of the conjugated operator. This term is compensated by the positive term $\Delta t\,|D_h^{2}D_tz|^{2}$ generated by the fourth-order part, at the price of a remainder $\Delta t\,s^{4}|D_tz|^{2}$, which is absorbed by $s^{-1}|D_tz|^{2}$ (Lemma \ref{lem:Dt:absorb}). These terms, which involve both $D_h$ and $D_t$ and carry a factor $\Delta t$, come from the interaction between the discrete spatial operator and the time discretization scheme, and therefore appear neither in the continuous nor in the semi-discrete setting. The last term in \eqref{eq:intro:carleman}, of order $h^{-4}$, collects the contributions at the initial and final times, which have no continuous counterpart.\\

From Theorem \ref{thm:intro:carleman} we deduce a relaxed observability inequality for \eqref{sys:adjoint}. For $\mu\ge1$ we set
$$
m_\mu:=\min\Big\{1,\frac{\mu}{5}\Big\},\qquad \Lambda:=1+\frac1T+\gamma^{2/7},\qquad h_1:=C_0\,\Lambda^{-\max\{1,5/\mu\}},\qquad C_{\rm obs}:=e^{C(\Lambda+T\gamma)}.
$$

\begin{theorem}\label{thm:intro:observability}
Let $\mu\ge1$. There exist $h_0,C_0,C,C_1>0$, depending on $\mu$ but independent of $T$, $a$, $h$ and $\Delta t$, such that, if
\begin{equation}\label{eq:intro:dt}
h\le\min\{h_0,h_1\}\qquad\text{and}\qquad \Delta t\le\min\Big\{h^{\mu},\ \frac{1}{4\gamma},\ \frac T4\Big\},
\end{equation}
then every solution of \eqref{sys:adjoint} satisfies
\begin{equation}\label{eq:intro:obs}
\|\trp q(0)\|^{2}_{L^2_h(\mathcal{M})}\le C_{\rm obs}\bigg(\int_{\mathcal{N}^{\ast}}|D_h^{2}q|^{2}\Big|_{x=0}+\int_{\mathcal{N}^{\ast}}|D_h^{3}q|^{2}\Big|_{x=h/2}+e^{-C_1/h^{m_\mu}}\|q_T\|^{2}_{L^2_h(\mathcal{M})}\bigg).
\end{equation}
\end{theorem}

The penalized Hilbert Uniqueness Method \cite{boyer-canum} then yields our main controllability result.

\begin{theorem}\label{thm:intro:control}
Under the assumptions of Theorem \ref{thm:intro:observability}, set $\phi(h):=e^{-C_1/h^{m_\mu}}$. For every $y_0\in L^2_h(\mathcal{M})$ there exist controls $u_1,u_2$ defined on $\mathcal{N}^{\ast}$ such that
$$
\|u_1\|_{L^2_{\Delta t}(\mathcal{N}^{\ast})}^{2}+\|u_2\|_{L^2_{\Delta t}(\mathcal{N}^{\ast})}^{2}\le C_{\rm obs}\|y_0\|^{2}_{L^2_h(\mathcal{M})},
$$
and the solution of \eqref{sys:primal} satisfies
$$
\|y(T)\|_{L^2_h(\mathcal{M})}\le\sqrt{C_{\rm obs}}\,\sqrt{\phi(h)}\,\|y_0\|_{L^2_h(\mathcal{M})}.
$$
\end{theorem}

\begin{remark}\label{rem:intro:conditions}
The scheme \eqref{sys:primal} is unconditionally stable, so neither condition on $\Delta t$ in \eqref{eq:intro:dt} is a stability (CFL) condition. The condition $\Delta t\,\gamma\le\frac14$ is only used to solve each implicit step and in a discrete Gr\"onwall argument for the adjoint system, and $\Delta t\le h^{\mu}$ comes from the second condition in \eqref{eq:intro:conditions}. In the proof, $\tau$ is fixed in terms of $T$ and $\gamma$, and $\delta$ is chosen proportional to $h^{m_\mu}$, so that $\tau(\delta T^{2})^{-1}$ is of order $h^{-m_\mu}$; this produces the exponentially small factor $\phi(h)$. The parameter $\mu$ quantifies how small $\Delta t$ is compared with $h$. Indeed, the rate $\phi(h)=e^{-C_1/h}$ is reached for $\mu\ge5$, whereas for $1\le\mu<5$ one obtains the slower rate $e^{-C_1/h^{\mu/5}}$. As in \cite{boyer-canum,AAP-2024}, Theorem \ref{thm:intro:control} also holds for any function $\phi$ satisfying $\liminf_{h\to0}\phi(h)\,e^{C_1/h^{m_\mu}}>0$.
\end{remark}

\par \textbf{Organization of the paper.} Section \ref{sec:prelim} collects the discrete calculus formulas in space and time, the duality identity between \eqref{sys:primal} and \eqref{sys:adjoint}, and the properties of the weight functions. Section \ref{sec:proof:carleman} is devoted to the proof of Theorem \ref{thm:intro:carleman}. In Section \ref{sec:control} we prove Theorems \ref{thm:intro:observability} and \ref{thm:intro:control}. Section \ref{sec:conclusion} compares our approach with \cite{kumar:2026,LW:2026} and lists some open problems. The estimates of the cross-product terms involving the discrete time derivative are given in Appendix \ref{estimaciones:carleman}.

\section{Preliminary results}\label{sec:prelim}
\subsection{Discrete calculus in space}\label{subsec:space}

We follow the notation of \cite[Section 1.1]{LOP:2021}. Let $\mathcal W$ be one of the spatial meshes introduced in Section \ref{subsec:setting}, that is, $\mathcal{M}$, $\overline{\mathcal M}$, $\widehat{\mathcal M}$, $\mathcal M^{\ast}$ or $\widehat{\mathcal M}^{\ast}$. We denote by $C(\mathcal W)$ the set of real functions defined on $\mathcal W$. For $u,v\in C(\mathcal W)$, the discrete integral, the inner product and the norm are
$$
\int_{\mathcal W}u:=h\sum_{x\in\mathcal W}u(x),\qquad \langle u,v\rangle_{L^2_h(\mathcal W)}:=\int_{\mathcal W}u\,v,\qquad \|u\|_{L^2_h(\mathcal W)}^2:=\int_{\mathcal W}|u|^2 .
$$
On the primal mesh we simply write $\langle u,v\rangle:=\langle u,v\rangle_{L^2_h(\mathcal{M})}$ and $\|u\|:=\|u\|_{L^2_h(\mathcal{M})}$. We use the translations $\tau_\pm(x):=x\pm\frac h2$ and $\tau_\pm u(x):=u(x\pm\frac h2)$. Throughout, $\mathcal W$ denotes one of the meshes $\mathcal{M}$ or $\overline{\mathcal M}$. Its closure is $\overline{\mathcal W}$, that is, $\overline{\mathcal M}$ or $\widehat{\mathcal M}$, its boundary is $\partial\mathcal W:=\overline{\mathcal W}\setminus\mathcal W$, that is, $\partial\mathcal{M}$ or $\partial\overline{\mathcal M}$, and its dual mesh is $\mathcal W^\ast:=\tau_+(\mathcal W)\cup\tau_-(\mathcal W)$, that is, $\mathcal M^{\ast}$ or $\widehat{\mathcal M}^{\ast}$. The boundary carries no factor $h$, namely
$$
\int_{\partial\mathcal W}u:=\sum_{x\in\partial\mathcal W}u(x).
$$
As in \cite[Section 1.1]{LOP:2021} and \cite[Section 1.2.1]{LMPZ:2023}, the outward normal at $x\in\partial\mathcal W$ is defined by
$$
n(x):=\begin{cases}
1, & \tau_-(x)\in\mathcal W^\ast\text{ and }\tau_+(x)\notin\mathcal W^\ast,\\
-1, & \tau_-(x)\notin\mathcal W^\ast\text{ and }\tau_+(x)\in\mathcal W^\ast,
\end{cases}
$$
so that $n=-1$ at the left boundary node and $n=1$ at the right one. For $g\in C(\mathcal W^\ast)$, the trace operator is defined, for $x\in\partial\mathcal W$, by
$$
t_r(g)(x):=\begin{cases}
\tau_-g(x), & n(x)=1,\\
\tau_+g(x), & n(x)=-1,
\end{cases}
$$
that is, $t_r(g)$ is the value of $g$ at the node of $\mathcal W^\ast$ adjacent to the boundary node.

We define the average operator $A_h$ and the difference operator $D_h$ as in Section \ref{subsec:setting}. The following product rules hold, with $\mathcal W$, $\overline{\mathcal W}$ and $\mathcal W^\ast$ as above.
\begin{lemma}[{\cite[Lemmas 3.2 and 3.3]{BHL10a}, \cite[Lemma 2.1]{LOP:2021}, \cite[Lemma 1]{LMPZ:2023}}]\label{pro:derivadapromedio}
Let $\mathcal W\in\{\mathcal{M},\overline{\mathcal M}\}$. For $u,v\in C(\overline{\mathcal W})$, the difference and average operators satisfy
\begin{equation}\label{eq:difference:product}
 D_{h}(u\,v)=D_{h}u\,A_{h}v+A_{h}u\,D_{h}v,\qquad
 A_{h}(u\,v)=A_{h}u\,A_{h}v+\frac{h^2}{4}D_{h}u\,D_{h}v,\qquad\text{in }\mathcal W^\ast,
\end{equation}
and, for $u\in C(\overline{\mathcal W})$,
\begin{equation}\label{eq:average:square}
u=A_{h}^{2}u-\frac{h^2}{4}D_{h}^{2}u\qquad\text{in }\mathcal W.
\end{equation}
The identities \eqref{eq:difference:product} also hold, with the same proof, for functions defined on a dual mesh, at the nodes of the corresponding primal mesh where both sides are defined.
\end{lemma}

As a direct consequence of Lemma \ref{pro:derivadapromedio} we have the following result.
\begin{corollary}[{\cite[Corollary 2.2]{CLTP:2022}}, {\cite[Corollary 2.2]{LOP:2021}}]\label{cor:identities}
Let $\mathcal W\in\{\mathcal{M},\overline{\mathcal M}\}$. For $u\in C(\overline{\mathcal W})$, we have, in $\mathcal W^\ast$,
\begin{equation}\label{promedio}
A_{h}(u^{2})=\left(A_{h}u \right)^{2}+\frac{h^{2}}{4}\left( D_{h}u\right)^{2}\geq \left(A_{h}u \right)^{2},\qquad D_{h}(u^{2})=2D_{h}u\,A_{h}u.
\end{equation}
\end{corollary}

The discrete Leibniz rule for $D_h^n$ \cite[Proposition 2.2]{AAP-2024} reads
\begin{equation}\label{eq:leibniz}
D_h^n(uv)=\sum_{k=0}^{n}\binom nk D_h^{n-k}A_h^{k}u\;A_h^{n-k}D_h^{k}v .
\end{equation}

\begin{lemma}[Discrete integration by parts, {\cite[Proposition 2.4]{LOP:2021}}]\label{lem:ipp}
Let $\mathcal W\in\{\mathcal{M},\overline{\mathcal M}\}$. For $u\in C(\overline{\mathcal W})$ and $v\in C(\mathcal W^\ast)$, we have
\begin{equation}\label{eq:ipp}
\int_{\mathcal W}u\,D_hv=-\int_{\mathcal W^\ast}v\,D_hu+\int_{\partial\mathcal W}u\,t_r(v)\,n,\qquad
\int_{\mathcal W}u\,A_hv=\int_{\mathcal W^\ast}v\,A_hu-\frac h2\int_{\partial\mathcal W}u\,t_r(v).
\end{equation}
\end{lemma}
For $\mathcal W=\mathcal{M}$, \eqref{eq:ipp} is also \cite[Proposition 1]{LMPZ:2023}. For $\mathcal W=\overline{\mathcal M}$ it will be used with $v=0$ on $\partial\widehat{\mathcal M}^{\ast}$; then both boundary terms vanish and the integrals over $\widehat{\mathcal M}^{\ast}$ reduce to integrals over $\mathcal M^{\ast}$.

\begin{remark}\label{rmk:green}
For $u,v$  defined on $\widehat{\mathcal M}$, with $u=v=0$ on $\partial\mathcal{M}$ and $D_hu=D_hv=0$ on $\partial\widehat{\mathcal M}^{\ast}$, we have
\begin{equation}\label{eq:green}
\int_{\mathcal{M}}(D_h^4u)\,v=\int_{\overline{\mathcal M}}D_h^2u\,D_h^2v .
\end{equation}
In particular, $\displaystyle\int_{\mathcal{M}}(D_h^4u)\,u\ge0$.
\end{remark}
Using the previous discrete calculus tools, we prove the duality identity between \eqref{sys:primal} and \eqref{sys:adjoint}.
\begin{lemma}\label{lem:duality}
Let $y$ be a solution of \eqref{sys:primal} and $q$ a solution of \eqref{sys:adjoint}. Then
\begin{equation}\label{eq:duality}
\langle y(T),q_T\rangle-\langle y_0,\trp q(0)\rangle=\int_{\mathcal{N}^{\ast}}\Big(u_2\,D_h^2q\big|_{x=0}-u_1\,D_h^3q\big|_{x=h/2}\Big).
\end{equation}
\end{lemma}
\begin{proof}
Fix $t\in\mathcal{N}^{\ast}$, notice that $\trp y(1)=0$, $D_h\trp y(1+\frac h2)=0$, $\trp y(0)=u_1(t)$, $D_h\trp y(-\frac h2)=u_2(t)$. Since $q(t)$ satisfies the homogeneous boundary conditions, repeated application of the first identity from \eqref{eq:ipp} gives $$\displaystyle\int_{\mathcal{M}}q\,D_h^4\trp y=\int_{\overline{\mathcal M}}D_h^2\trp y\,D_h^2q.$$
On the other hand, $\trp y(-h)=\trp y(0)-hu_2$ yields $t_r(D_h\trp y)(0)=hD_h^2\trp y(0)+u_2$ and $t_r(D_h\trp y)(1)=-hD_h^2\trp y(1)$; applying \eqref{eq:ipp} twice in the reverse direction we obtain
$$
\int_{\mathcal{M}}q\,D_h^4\trp y=\int_{\mathcal{M}}\trp y\,D_h^4q+u_1(t)\,D_h^3q(\tfrac h2,t)-u_2(t)\,D_h^2q(0,t).
$$
Multiplying the first equation of \eqref{sys:primal} by $q$ and integrating over $\mathcal{M}\times\mathcal{N}^{\ast}$, we use \eqref{eq:tau:mas} for the time derivative, \eqref{eq:tau:menos} for the terms $\trp(\cdot)$, the identity above, and the first equation of \eqref{sys:adjoint}, and we obtain
$$
0=\langle y(T),q_T\rangle-\langle y_0,\trp q(0)\rangle+\int_{\mathcal{M}\times\mathcal{N}}y\big(-D_tq+D_h^4\trm q+a\,\trm q\big)+\int_{\mathcal{N}^{\ast}}\Big(u_1D_h^3q\big|_{x=h/2}-u_2D_h^2q\big|_{x=0}\Big),
$$
and the third term vanishes.
\end{proof}

\subsection{Time discrete calculus formulas}
In this section, we introduce the notation used in discrete-time calculus. The time meshes $\mathcal{N}$, $\overline{\mathcal N}$, $\mathcal{N}^{\ast}$, $\overline{\mathcal N}^{\ast}$, the translations $\mathtt t^\pm$ and $D_t$ were defined in Section \ref{subsec:setting}. The discrete integrals are $\displaystyle\int_{\mathcal{N}}f:=\Delta t\sum_{t\in\mathcal{N}}f(t)$, $\displaystyle\int_{\mathcal{N}^{\ast}}f:=\Delta t\sum_{t\in\mathcal{N}^{\ast}}f(t)$ and $\displaystyle\int_{\partial\mathcal{N}}f\,n:=f(T)-f(0)$, and $\displaystyle\|u\|^2_{L^2_{\Delta t}(\mathcal{N}^{\ast})}:=\int_{\mathcal{N}^{\ast}}|u|^2$. Similarly to the discrete spatial variable, we define the outward normal for $t\in \partial \mathcal{N}:=\{0,T\}$ as $n(T)=1$, $n(0)=-1$.

In what follows, we present some identities that will be useful in our main results; they are taken from \cite[Section 2.2]{LMPZ:2023}. For $f$ and $g$ continuously defined over $\mathbb{R}$, the following identities hold
\begin{equation}\label{eq:Dt:product}
    D_{t}(fg)=\trp f\,D_{t}g+D_{t}f\,\trm g=D_{t}f\,\trp g+D_{t}g\,\trm f,
\end{equation}
see \cite[(2.1)--(2.2)]{LMPZ:2023}. Moreover, considering $f=g$ we obtain \cite[(2.3)]{LMPZ:2023}
\begin{equation}\label{eq:identity:square}
       \trp f\,D_{t}f=\frac{1}{2}D_{t}(f^{2})+\frac{1}{2}\Delta t|D_{t}f|^{2},\qquad
        \trm f\,D_{t}f=\frac{1}{2}D_{t}(f^{2})-\frac{1}{2}\Delta t|D_{t}f|^{2}.
\end{equation}
Let us finally comment that we will use the above identities sampled on the primal and dual time-discrete meshes. Now, we present two integration by parts formulas for the time-discrete operator.
\begin{proposition}[{\cite[Proposition 2]{LMPZ:2023}}]
    \label{properties:time:discret:ops}
    Let us consider $f\in C(\overline{\mathcal N})$ and $g\in C(\overline{\mathcal N}^{\ast})$. Then, the following identities hold
    \begin{equation}\label{eq:tau:menos} \int_{\mathcal{N}} f\,\trm(g)=\int_{\mathcal{N}^{\ast}}\trp (f)\, g,
    \end{equation}
    \begin{equation}\label{eq:tau:mas}     \int_{\mathcal{N}}f\,D_{t}g=-\int_{\mathcal{N}^{\ast}}g\,D_{t}f+\int_{\partial \mathcal{N}}f\,\trp g\,n.
    \end{equation}
\end{proposition}
Moreover, by \cite[(2.7)]{LMPZ:2023},
\begin{equation}\label{time:integral}
\int_{\mathcal{N}}\trp g=\int_{\mathcal{N}^{\ast}}g+\Delta t\int_{\partial\mathcal{N}}\trp g\,n,
\end{equation}
and, for any $F,G$ defined on $\overline{\mathcal N}^{\ast}$, the time-slicing identity
\begin{equation}\label{eq:slicing}
\int_{\mathcal{N}}\trm F\,\trm G=\int_{\mathcal{N}^{\ast}}F\,G
\end{equation}
holds, since $\trm$ maps $\mathcal{N}$ onto $\mathcal{N}^{\ast}$; see \cite[(2.6)]{LMPZ:2023}.

\subsection{Weight functions}\label{subsec:weights}

As in the semi-discrete setting \cite[Equations (7)--(8)]{CLTP:2022}, the spatial part of the weight satisfies the following assumption.
\begin{assumption}\label{ass:psi}
Let $\widetilde I$ be an open interval containing $[-1,2]$. The function $\varphi\in C^k(\widetilde I)$, with $k$ large enough, satisfies, for some $\eta>0$,
$$
0<\eta\le\varphi(x),\qquad 0<\eta\le\frac{d\varphi}{dx}(x),\qquad \frac{d^2\varphi}{dx^2}(x)\le-\eta<0,\qquad \text{for all }x\in[-1,2].
$$
\end{assumption}
We set
$$
k_0:=\min_{[-1,2]}\varphi>0,\qquad K_0:=\max_{[-1,2]}\varphi,\qquad S:=\tau(\delta T^2)^{-1}.
$$
As in \cite[Section 2.2]{CLTP:2022}, $\mathcal{O}(1)$ denotes a function bounded by a constant that depends only on $\varphi$, and not on $s,\tau,h,\Delta t,T$ or $\delta$, and $\mathcal{O}(X):=X\mathcal{O}(1)$. All constants below may depend on $\varphi$.

We now collect the one-dimensional estimates of \cite[Section 2.2]{CLTP:2022} on the discrete operators applied to the weights; see also \cite[Section 2.3]{LMPZ:2023} and, in arbitrary dimension, \cite{AAP-2024}.

\begin{proposition}[{\cite[Lemma 2.5 and Theorem 2.8]{CLTP:2022}}]\label{prop:weights}
Let $k,l,m,n,\alpha\in\mathbb N_0$. Assume $\tau\ge\tau_0(T+T^2)$, $\delta\in(0,\frac12]$ and \eqref{eq:intro:conditions}. Then, at every point of $[-1,2]\times[-\Delta t,T+\Delta t]$ where they are defined,
\begin{align}
&A_h^lD_h^k\partial_x^\alpha\big(rA_h^mD_h^n\rho\big)=\partial_x^{k+\alpha}\big(r\partial_x^n\rho\big)+s^{n}\mathcal{O}((sh)^2),\label{eq:AAP:space}\\
&r\partial_x^n\rho=(s\,\partial_x\varphi)^n+s^{n-1}\mathcal{O}(1),\qquad \partial_x^k\big(r\partial_x^n\rho\big)=s^n\mathcal{O}(1).\label{eq:AAP:cont}
\end{align}
In particular $A_h^lD_h^k(rA_h^mD_h^n\rho)=s^n\mathcal{O}(1)$; $A_h^l(rA_h^m\rho)=1+\mathcal{O}((sh)^2)$; and $A_h^lD_h^k(rA_h^m\rho)=\mathcal{O}((sh)^2)$ for $k\ge1$.
\end{proposition}
\begin{remark}
Since $\tau\ge T$ and $\tau\ge T^2$, we have $\tau^5\ge T^9$, and the second condition in \eqref{eq:intro:conditions} gives
\begin{equation}\label{eq:dt:small}
\Delta t\le\varepsilon_0\,\delta^5T^{10}\tau^{-5}\le\varepsilon_0\,\delta^5T\le\frac{\delta T}{16}.
\end{equation}
Hence, for $-\Delta t\le t\le0$, $t+\delta T\ge\frac{15}{16}\delta T$ and $T+\delta T-t\ge T$, so that $\theta(t)\le2(\delta T^2)^{-1}$; the same bound holds for $T\le t\le T+\Delta t$ by symmetry, and on $[0,T]$ it holds since $\theta\le\theta(0)$. With $S=\tau(\delta T^2)^{-1}$ and the first condition in \eqref{eq:intro:conditions}, we obtain
\begin{equation}\label{eq:s:bound}
s\le2S\quad\text{and}\quad sh\le2\varepsilon_0\le1\qquad\text{on }[-\Delta t,T+\Delta t].
\end{equation}
Moreover $|\theta'|=|2t-T|\theta^2\le2T\theta^2$ there, so $|(\log\theta)'|\le2T\theta\le4(\delta T)^{-1}$, and by \eqref{eq:dt:small},
\begin{equation}\label{eq:theta:shift}
e^{-1/4}\theta(t)\le\theta(t')\le e^{1/4}\theta(t)\qquad\text{for }t,t'\in[-\Delta t,T+\Delta t],\ |t-t'|\le\Delta t .
\end{equation}

\end{remark}

\begin{lemma}[{\cite[Lemma 5]{LMPZ:2023}}, {\cite[Lemma B.4]{boyer-hernandez}}]\label{lem:time:conj}
Let $\mathcal E:=\Delta t\big(\tau(\delta^3T^4)^{-1}+\tau^2(\delta^4T^6)^{-1}\big)$. Assume $\tau\ge\tau_0(T+T^2)$, $\delta\in(0,\frac12]$ and \eqref{eq:intro:conditions}. Then, on $\mathcal{N}$,
$$
\trm(r)\,D_t\rho=\tau\,\varphi\,\trm(\theta')+\mathcal E\,\mathcal{O}(1),
$$
and $\mathcal E\le2\varepsilon_0$.
\end{lemma}

\begin{lemma}\label{lem:Dt:outer}
Let $k,l,m,n\in\mathbb N_0$, $F:=A_h^lD_h^k\big(rA_h^mD_h^n\rho\big)$ and $\beta\in\{0,\frac12,1\}$. Assume $\tau\ge\tau_0(T+T^2)$, $\delta\in(0,\frac12]$ and \eqref{eq:intro:conditions}. Then, for every $t$ such that $t-\beta\Delta t,\,t+(1-\beta)\Delta t\in[-\Delta t,T+\Delta t]$,
$$
\frac{F(t+(1-\beta)\Delta t)-F(t-\beta\Delta t)}{\Delta t}=T\,\theta(t)\,s(t)^{n}\,\mathcal{O}(1).
$$
In particular $D_tF$ and $D_t(\trm F)$ are $T\theta s^n\mathcal{O}(1)$.
\end{lemma}
This is the time-discrete version of \cite[Theorem 2.9, (29)]{CLTP:2022}; for $k=l=0$ it is also \cite[Theorem 5]{LMPZ:2023}.
\begin{proof}
The difference quotient equals $\displaystyle\int_0^1\partial_tF(t')\,d\gamma$ with $t'=t-\beta\Delta t+\gamma\Delta t$. By \cite[Theorem 2.9, (29)]{CLTP:2022}, written with the roles of $r$ and $\rho$ interchanged, $\partial_tF=T\theta s^n\mathcal{O}(1)$ at every time $t'$ such that $s(t')h\le1$, and by \eqref{eq:s:bound} this holds on $[-\Delta t,T+\Delta t]$. Since $|t'-t|\le\Delta t$, \eqref{eq:theta:shift} gives $\theta(t')\le e^{1/4}\theta(t)$ and $s(t')\le e^{1/4}s(t)$, which proves the claim.
\end{proof}

\section{Proof of Theorem \ref{thm:intro:carleman}}\label{sec:proof:carleman}

Throughout this section we assume $\tau\ge\tau_0(T+T^2)$ and \eqref{eq:intro:conditions}, with $\tau_0$ large and $\varepsilon_0$ small to be fixed at the end; hence Proposition \ref{prop:weights} and Lemmas \ref{lem:time:conj} and \ref{lem:Dt:outer} apply, and \eqref{eq:dt:small}--\eqref{eq:theta:shift} hold. $C$ denotes a generic positive constant independent of $\tau,h,\Delta t,T,\delta,\tau_0,\varepsilon_0$. We follow the traditional approach of conjugating the original operator with a carefully chosen exponential function. We write the conjugate operator into a symmetric part $S_{h}$, an antisymmetric part $P_{h}$, and an additional term $R_{h}$. Then, we estimate the cross-inner product between these operators, handling the purely spatial terms time slice by time slice as in \cite{CLTP:2022}, and the terms involving $D_t$ in Appendix \ref{estimaciones:carleman}; as a final stage, we return to the original variable.

\subsection{Conjugation with the weight function}
Let $q$ be defined on $\widehat{\mathcal M}\times\overline{\mathcal N}^{\ast}$ and satisfy the boundary conditions of \eqref{sys:adjoint} at every $t\in\overline{\mathcal N}^{\ast}$. We set $z:=rq$, that is, $q=\rho z$. Then $z=0$ on $(\partial\mathcal{M}\cup\partial\overline{\mathcal M})\times\overline{\mathcal N}^{\ast}$, and $D_{t}z=0$ on $(\partial\mathcal{M}\cup\partial\overline{\mathcal M})\times\mathcal{N}$; in particular $D_hz=D_hD_{t}z=A_hD_{t}z=0$ on $\partial\widehat{\mathcal M}^{\ast}$.
Applying \eqref{eq:leibniz} it follows that
\begin{equation}
\label{eq:fourth:dh}
\begin{aligned}D_{h}^4(\rho z)=&D_{h}^{4}(\rho)\, A_{h}^{4}(z)+4D_{h}^{3}A_{h}(\rho)\, D_{h} A_{h}^{3}(z)+6A_{h}^{2}D_{h}^{2}(\rho)\, D_{h}^{2} A_{h}^{2}(z)\\
&+4 A_{h}^{3}D_{h}(\rho)\, D_{h}^{3}A_{h}(z)+A_{h}^{4}(\rho)\,D_{h}^{4}(z).
\end{aligned}
\end{equation}
For the time derivative, \eqref{eq:Dt:product} gives $D_t(\rho z)=D_t\rho\,\trp z+\trm\rho\,D_tz$; by Lemma \ref{lem:time:conj} and $\trp z=\trm z+\Delta t\,D_tz$,
\begin{equation}\label{eq:change_discr_temp}
-\trm(r)D_{t}(\rho z)=-D_tz-\tau\varphi\,\trm(\theta' z)-\Delta t\,\tau\varphi\,\trm(\theta')D_tz+\mathcal E\,\mathcal{O}(1)\trp z .
\end{equation}
Thus, from the change of variable $q=\rho z$, with $\mathcal P^\star q=-D_tq+D_h^4\trm q$, we obtain
$$
S_h z + P_h z = \trm r\, \mathcal{P}^{\star} (\rho z) + R_h z ,
$$
where $S_hz=\sum_{i=1}^4 S_i z$, $P_hz=\sum_{i=1}^{4}{P_i} z$,
\begin{align*}
&S_1 z= 6\trm\left(rA_{h}^{2}D_{h}^{2}(\rho)\, D_{h}^{2} A_{h}^{2}(z)\right), \quad S_2 z=\trm\left(rD_{h}^{4}(\rho)\, A_{h}^{4}(z)\right), \quad S_3 z= \trm\left(rA_{h}^{4}(\rho)\, D_{h}^{4}(z)\right),\\
& S_4 z=6\trm\left(A_{h}D_{h}(r A_{h}^2D_{h}^2\rho)\,A_{h}D_{h}z\right),\\
&P_1z=-D_{t} z, \quad P_2 z=4\trm\left(rD_{h}^{3}A_{h}(\rho)\, D_{h} A_{h}^{3}(z)\right), \quad P_3z=4 \trm\left(rA_{h}^{3}D_{h}(\rho)\, D_{h}^{3}A_{h}(z)\right), \\
&P_4 z = 2\trm\left(A_{h}D_{h}(r \,A_{h}D_{h}^3\rho)\,z\right),
\end{align*}
and
\begin{equation}\label{def:R_term}
R_hz:=\tau\varphi\,\trm(\theta'z)+\Delta t\,\tau\varphi\,\trm(\theta')D_tz-\mathcal E\,\mathcal{O}(1)\trp z+S_4z+P_4z .
\end{equation}
Then, 
\begin{equation}\label{eq:general estimate}
\left\|{S_{h}z}\right\|^{2}_{L_{h}^{2}(Q)}+\left\|{P_{h}z} \right\|^{2}_{L_{h}^{2}(Q)}+2\int_{Q} S_{h}z\, P_{h}z
=\left\| {\trm r\, \mathcal{P}^{\star}q+R_{h}z}\right\|^{2}_{L_{h}^{2}(Q)}
\le2\left\|{\trm r\, \mathcal{P}^{\star}q} \right\|^{2}_{L_{h}^{2}(Q)} +2\left\|{R_{h}z} \right\|^{2}_{L_{h}^{2}(Q)}.
\end{equation}
We will write
\begin{equation}\label{eq:Qmeshes}
Q:= \mathcal{M} \times \mathcal{N},\quad Q^{\ast}:= \mathcal M^{\ast} \times \mathcal{N}^{\ast},\quad Q_{\ast}:= \mathcal{M} \times \mathcal{N}^{\ast},\quad \overline{Q}_\ast := \overline{\mathcal M} \times \mathcal{N}^{\ast} .
\end{equation}
We also introduce the quantities
\begin{equation}\label{eq:quantities}
\begin{aligned}
&L(z):=\int_{Q_{\ast}}s^{7}|z|^{2}+\int_{Q^{\ast}}s^{5}|D_{h}z|^{2}+\int_{\overline{Q}_{\ast}}s^{3}|D_{h}^{2}z|^{2}+\int_{Q^{\ast}}s|D_{h}^{3}z|^{2},\\
&\mathcal{D}(z):=\int_{Q}\trm(s^{-1})|D_tz|^2,\qquad \mathcal{F}(z):=\int_{Q_{\ast}}s^{-1}|D_h^4z|^2,\qquad \mathcal{T}(z):=\int_{\mathcal{M}\times\partial\mathcal{N}}|\trp z|^2,\\
&\mathcal{B}(z):=\int_{\mathcal{N}^{\ast}}\Big(s^3|D_h^2z|^2\big|_{x=0}+s|D_h^3z|^2\big|_{x=h/2}+s^5|D_hz|^2\big|_{x=h/2}\Big).
\end{aligned}
\end{equation}

\subsection{An estimate for the remainder}

\begin{lemma}\label{lem:estimate Rv}
We have
$$
\left\|{R_{h}z} \right\|^2_{L_{h}^{2}(Q)} \le C\tau_0^{-1}L(z)+C\varepsilon_0^2\mathcal{D}(z)+Ch^{-4}\mathcal T(z).
$$
\end{lemma}
\begin{proof}
By Proposition \ref{prop:weights}, $A_hD_h(rA_h^2D_h^2\rho)=s^2\mathcal{O}(1)$ and $A_hD_h(rA_hD_h^3\rho)=s^3\mathcal{O}(1)$; using \eqref{eq:slicing}, $|A_hD_hz|^2\le A_h(|D_hz|^2)$ and \eqref{eq:ipp},
$$
\|S_4z\|^2+\|P_4z\|^2\le C\int_{Q^{\ast}}s^4|D_hz|^2+C\int_{Q_{\ast}}s^6|z|^2\le C\tau_0^{-1}L(z).
$$
Since $|\theta'|\le2T\theta^2$ and $\tau\ge T$, we have $T\theta=Ts/\tau\le s$ and $|\tau\varphi\theta'|\le CTs\theta\le Cs^2$; moreover $\theta\ge4(1+2\delta)^{-2}T^{-2}\ge T^{-2}$ on $[0,T]$, so $s\ge\tau T^{-2}\ge\tau_0$. Hence the first term of $R_h$ contributes $\displaystyle C\int_{Q_{\ast}}s^4|z|^2\le C\tau_0^{-1}L(z)$. For the second term, $\Delta t^2\tau^2|\theta'|^2\le\Delta t^2\tau^3T^2\theta^5\,s^{-1}\le\Delta t^2\tau^3(\delta^5T^8)^{-1}s^{-1}$, and by \eqref{eq:intro:conditions} and $\tau^7\ge T^{12}$, $\Delta t^2\tau^3(\delta^5T^8)^{-1}\le\varepsilon_0^2\delta^5T^{12}\tau^{-7}\le\varepsilon_0^2$. This gives $C\varepsilon_{0}^{2}\mathcal{D}(z)$. Finally, by \eqref{time:integral}, $\displaystyle\int_Q|\trp z|^2\le\int_{Q_{\ast}}|z|^2+\Delta t\,\mathcal{T}(z)$, and we use $\mathcal E\le2\varepsilon_0$ (Lemma \ref{lem:time:conj}), $s\ge\tau_0$, and $\Delta t\le\varepsilon_0\delta^5T^{10}\tau^{-5}\le1\le h^{-4}$, which follows from \eqref{eq:intro:conditions} and $\tau^5\ge T^{10}$.
\end{proof}

\subsection{An estimate for the cross term}
In this stage, we provide an estimate for the terms
\begin{equation}\label{eq:conmutador}
    \int_{Q}S_{h}z\,P_{h}z=\sum_{i=1}^{4}\sum_{j=1}^{4}I_{ij},
\end{equation}
where $I_{ij}$ stands for the inner product in $L^2_{h}(Q)$ between the $i^{\text{th}}$-term of $S_{h}z$ and the $j^{\text{th}}$-term of $P_{h} z$.

\begin{proposition}[Spatial cross terms]\label{prop:spatial}
There exist $c_0>0$ and $C>0$ such that
$$
\sum_{i=1}^{4}\sum_{j=2}^{4}I_{ij}\ge c_0L(z)-C\big(\varepsilon_0+\tau_0^{-1}\big)L(z)-CB(z).
$$
\end{proposition}
\begin{proof}
For $j\ge2$, $S_iz$ and $P_jz$ are of the form $\trm(\cdot)$ of purely spatial expressions. Hence, by \eqref{eq:slicing}, $\displaystyle I_{ij}=\int_{\mathcal{N}^{\ast}}\mathcal I_{ij}(t)$, where $\mathcal I_{ij}(t)$ is the corresponding spatial integral over $\mathcal{M}$ at the fixed time $t\in\mathcal{N}^{\ast}$, with the weights $r(\cdot,t)=e^{-s(t)\varphi}$ and $\rho(\cdot,t)=e^{s(t)\varphi}$. With the roles of $r$ and $\rho$ interchanged, the operators $S_i$ and $P_j$ coincide with those of \cite[(32)]{CLTP:2022}, the function $z(\cdot,t)$ satisfies the boundary conditions of the semi-discrete adjoint system \cite[(6)]{CLTP:2022}, and $\mathcal I_{ij}(t)$ is the integrand in time of the term $I_{ij}$ of \cite[(40)]{CLTP:2022}. These terms are estimated in \cite[Lemmas 6.2--6.4, 6.6--6.8 and 6.10--6.15]{CLTP:2022}, for $I_{12},I_{13},I_{14}$, $I_{22},I_{23},I_{24}$, $I_{32},I_{33},I_{34}$ and $I_{42},I_{43},I_{44}$, respectively. Their proofs only use Lemmas \ref{pro:derivadapromedio} and \ref{lem:ipp} in space and the estimates of Proposition \ref{prop:weights}, at a fixed time, under the condition $s(t)h\le1$; no integration by parts in time is involved, since $P_1=-D_t$ only enters the terms $I_{i1}$. By \eqref{eq:s:bound}, $s(t)h\le2\varepsilon_0$ for every $t\in\mathcal{N}^{\ast}$, so these estimates hold at each $t\in\mathcal{N}^{\ast}$, and summing over $\mathcal{N}^{\ast}$ we obtain, as in \cite[(41)--(46)]{CLTP:2022},
$$
\sum_{i=1}^{4}\sum_{j=2}^{4}I_{ij}\ge\mathcal L(z)+\sum_{i=1}^4\sum_{j=2}^4X_{ij}+Y,
$$
where $\mathcal L(z)$ collects the leading terms of Remark \ref{rem:leading} together with terms of lower order in $s$, namely
$$
\begin{aligned}
&8\int_{Q_{\ast}}s^5\partial_x^2\big((\partial_x\varphi)^4\partial_x^2\varphi\big)|z|^2+3\int_{Q_{\ast}}s^3\partial_x^4\big((\partial_x\varphi)^2\partial_x^2\varphi\big)|z|^2-36\int_{Q_{\ast}}s^5\partial_x\big((\partial_x\varphi)^3(\partial_x^2\varphi)^2)\big)|z|^2\\
&\qquad+6\int_{Q^{\ast}}s^3\partial_x^2\big((\partial_x\varphi)^2\partial_x^2\varphi\big)|D_hz|^2,
\end{aligned}
$$
the terms $X_{ij}$ are the remainders of \cite[Lemmas 6.2--6.15]{CLTP:2022}, and $Y$ is the sum of the spatial boundary terms $Y_{12},Y_{13},Y_{22},Y_{23},Y_{32},Y_{33},Y_{42},Y_{43}$ of \cite[(46)]{CLTP:2022}. Since $s\ge\tau_0$ on $[0,T]$, the lower-order terms of $\mathcal L(z)$ and every $X_{ij}$, whose coefficients are of the form $s^{m}\mathcal{O}(1)$ or $s^{m+1}\mathcal{O}((sh)^2)$ in front of a term of $L(z)$ with weight $s^{m+1}$, are bounded by $C(\varepsilon_0^2+\tau_0^{-1})L(z)$. Finally, as shown in \cite[(48)--(49)]{CLTP:2022}, after summation the leading parts of the boundary terms in $Y$ are $s^5(\partial_x\varphi)^5t_r(|D_hz|^2)n$, $s^3(\partial_x\varphi)^3|D_h^2z|^2n$ and $s\,\partial_x\varphi\,t_r(|D_h^3z|^2)n$, multiplied by positive constants, and the remaining parts carry an additional factor $\mathcal{O}((sh)^2)$ or a lower power of $s$. At $x=1$ we have $n=1$ and $\partial_x\varphi>0$, so these terms are nonnegative for $\varepsilon_0$ small and $\tau_0$ large, and they can be dropped. At $x=0$ we have $n=-1$, and these terms are bounded by $C\mathcal{B}(z)$, since $t_r(|D_hz|^2)(0)=|D_hz(\frac h2)|^2$ and $t_r(|D_h^3z|^2)(0)=|D_h^3z(\frac h2)|^2$. Together with Remark \ref{rem:leading}, this proves the claim.
\end{proof}

\begin{remark}\label{rem:leading}
The leading terms produced by the spatial cross terms are, as in \cite[(42)--(45)]{CLTP:2022},
$$
\begin{aligned}
&-8\int_{Q_{\ast}}s^7(\partial_x\varphi)^6\partial_x^2\varphi\,|z|^2,\qquad -18\int_{Q^{\ast}}s^5(\partial_x\varphi)^4\partial_x^2\varphi\,|D_hz|^2,\\
&-60\int_{\overline{Q}_{\ast}}s^3(\partial_x\varphi)^2\partial_x^2\varphi\,|D_h^2z|^2,\qquad -2\int_{Q^{\ast}}s\,\partial_x^2\varphi\,|D_h^3z|^2 .
\end{aligned}
$$
They are positive and bounded below by $c_0L(z)$, with $c_0:=2\min\{\eta,\eta^7\}$, since $\partial_x\varphi\ge\eta$ and $-\partial_x^2\varphi\ge\eta$ by Assumption \ref{ass:psi}.
\end{remark}

\begin{proposition}[Cross terms with the time derivative]\label{prop:time}
There exists $C_1>0$ such that
$$
\sum_{i=1}^{4}I_{i1}\ge\frac{\Delta t}{4}\int_{\overline{\mathcal M}\times\mathcal{N}}|D_h^2D_tz|^2-C_1\Delta t\int_{\mathcal M^{\ast}\times\mathcal{N}}\trm(s^2)|D_hD_tz|^2-C\big(\varepsilon_0+\tau_0^{-1}\big)\big(L(z)+\mathcal{D}(z)\big)-Ch^{-4}\mathcal T(z).
$$
\end{proposition}
\begin{proof}
Add Lemmas \ref{lem:estimate:I11}, \ref{lem:estimate:I21} and \ref{lem:estimate:I31} of Appendix \ref{estimaciones:carleman}. Since $\tau\ge\tau_0T$, we have $T\theta=Ts/\tau\le s/\tau_0$, which allows us to bound $\displaystyle\int T\theta s^4|z|^2$, $\displaystyle\int T\theta s^2|D_hz|^2$ and $\displaystyle\int T\theta|D_h^2z|^2$ by $\tau_0^{-1}L(z)$.
\end{proof}

\subsection{Absorbing the time-discretization term}

\begin{lemma}\label{lem:Dt:absorb}
For every $C_1>0$,
$$
C_1\Delta t\int_{\mathcal M^{\ast}\times\mathcal{N}}\trm(s^2)|D_hD_tz|^2\le\frac{\Delta t}{8}\int_{\overline{\mathcal M}\times\mathcal{N}}|D_h^2D_tz|^2+64\,C_1^2\,\varepsilon_0\,\mathcal{D}(z).
$$
\end{lemma}
\begin{proof}
Fix $t\in\mathcal{N}$ and notice that $D_tz$ vanishes on $\partial\mathcal{M}\cup\partial\overline{\mathcal M}$. The first identity in \eqref{eq:ipp} and the Cauchy--Schwarz inequality, gives
$$
\int_{\mathcal M^{\ast}}|D_hD_{t}z|^2=-\int_{\mathcal{M}}D_{t}z\,D_h^2D_{t}z\le\|D_{t}z\|_{L^2_h(\mathcal{M})}\,\|D_h^2D_{t}z\|_{L^2_h(\overline{\mathcal M})} .
$$
Combining this bound with Young's inequality, we get
$$
C_1\trm(s^2)\int_{\mathcal M^{\ast}}|D_hD_{t}z|^2\le C_1\trm(s^2)\|D_{t}z\|\,\|D_h^2D_{t}z\|\le\frac18\|D_h^2D_{t}z\|^2_{L^2_h(\overline{\mathcal M})}+2C_1^2\trm(s^4)\|D_{t}z\|^2_{L^2_h(\mathcal{M})}.
$$
Multiply by $\Delta t$ and sum over $\mathcal{N}$. By \eqref{eq:s:bound} and the second condition in \eqref{eq:intro:conditions}, $\Delta t\,\trm(s^5)\le32\,\Delta t\,S^5=32\,\tau^5\Delta t(\delta^5T^{10})^{-1}\le32\varepsilon_0$, hence $\Delta t\,\trm(s^4)=\Delta t\,\trm(s^5)\trm(s^{-1})\le32\varepsilon_0\trm(s^{-1})$.
\end{proof}

\subsection{Recovering the time derivative and the fourth-order term}
\begin{lemma}\label{lem:absorb}
$\mathcal{D}(z)\le2\|P_hz\|^2_{L^2_h(Q)}+C\,L(z)$.
\end{lemma}
\begin{proof}
$D_tz=-P_hz+P_2z+P_3z+P_4z$ and $\trm(s^{-1})\le1$. By Proposition \ref{prop:weights}, the coefficients of $P_2$, $P_3$, $P_4$ are $s^3\mathcal{O}(1)$, $s\,\mathcal{O}(1)$, $s^3\mathcal{O}(1)$. Recalling that $|A_h^kf|^2\le A_h^k(|f|^2)$, and using the second identity in \eqref{eq:ipp} we have $\displaystyle\|\trm(s^{-1/2})(P_2+P_3+P_4)z\|^2\le C\int_{Q_{\ast}}\big(s^5|D_hA_h^3z|^2+s|D_h^3A_hz|^2+s^5|z|^2\big)\le CL(z)$, where we have used \eqref{eq:slicing} for the time translation operator.
\end{proof}

\begin{lemma}\label{lem:fourth:term}
$\mathcal{F}(z)\le8\|S_hz\|^2_{L^2_h(Q)}+C\,L(z)$.
\end{lemma}
\begin{proof}
By Proposition \ref{prop:weights}, $rA_h^4\rho=1+\mathcal{O}((sh)^2)\ge\frac12$, so $\mathcal{F}(z)\le4\|\trm(s^{-1/2})S_3z\|^2$. Since $S_3z=S_hz-S_1z-S_2z-S_4z$, with coefficients $s^2\mathcal{O}(1)$, $s^4\mathcal{O}(1)$, $s^2\mathcal{O}(1)$, and $\trm(s^{-1})\le1$,
$$
\|\trm(s^{-1/2})S_3z\|^2\le2\|S_hz\|^2+C\int_{Q_{\ast}}\big(s^3|D_h^2A_h^2z|^2+s^7|A_h^4z|^2+s^3|A_hD_hz|^2\big).
$$
By the bound $|A_h^kf|^2\le A_h^k(|f|^2)$, $\displaystyle\int_{\mathcal{M}}|D_h^2A_h^2z|^2=\int_{\mathcal{M}}|A_h^2D_h^2z|^2\le\int_{\overline{\mathcal M}}|D_h^2z|^2$, $\displaystyle\int_{\mathcal{M}}|A_h^4z|^2\le\int_{\mathcal{M}}|z|^2$ and $\displaystyle\int_{\mathcal{M}}|A_hD_hz|^2\le\int_{\mathcal M^{\ast}}|D_hz|^2$, since $z=0$ on $\partial\mathcal{M}\cup\partial\overline{\mathcal M}$ and $D_hz=0$ on $\partial\widehat{\mathcal M}^{\ast}$.
\end{proof}

\subsection{Boundary terms and return to the original variable}

\begin{lemma}\label{lem:boundary}
$\mathcal{B}(z)\le C\displaystyle\int_{\mathcal{N}^{\ast}}r^2\big|_{x=0}\Big(s^3|D_h^2q|^2\big|_{x=0}+s|D_h^3q|^2\big|_{x=h/2}\Big)$.
\end{lemma}
\begin{proof}
Fix $t\in\mathcal{N}^{\ast}$. Since $q(0)=q(-h)=0$ and $z(0)=z(-h)=0$,
$$
D_h^2z(0)=\frac{z(h)}{h^2}=r(h)D_h^2q(0),\qquad D_hz(\tfrac h2)=\frac{z(h)}h=h\,D_h^2z(0),\qquad D_h^3z(\tfrac h2)=\frac{z(2h)-3z(h)}{h^3}.
$$
Writing $q(2h)=h^3D_h^3q(\frac h2)+3h^2D_h^2q(0)$, we get
$$
D_h^3z(\tfrac h2)=r(h)D_h^3q(\tfrac h2)+\frac{r(2h)-r(h)}{h^3}\Big(h^3D_h^3q(\tfrac h2)+3h^2D_h^2q(0)\Big),
$$
and, $r(x')/r(x)=e^{-s(\varphi(x')-\varphi(x))}$ with $s|\varphi(x')-\varphi(x)|\le2sh\|\varphi'\|_\infty\le4\varepsilon_0\|\varphi'\|_\infty$ for $|x-x'|\le2h$, by \eqref{eq:s:bound}. Therefore $|r(2h)-r(h)|\le Csh\,r(h)$ and $r(h)\le Cr(0)$. Hence $$s|D_h^3z(\frac h2)|^2\le Cr(0)^2\big(s|D_h^3q(\frac h2)|^2+s^3|D_h^2q(0)|^2\big),$$ and  $$s^5|D_hz(\frac h2)|^2=(sh)^2s^3|D_h^2z(0)|^2\le Cs^3r(0)^2|D_h^2q(0)|^2.$$
\end{proof}

\begin{lemma}\label{lem:primalvariable}
The left-hand side of \eqref{eq:intro:carleman} is bounded by $C\big(L(z)+\mathcal{D}(z)+\mathcal{F}(z)+h^{-4}\mathcal{T}(z)\big)$.
\end{lemma}
\begin{proof}
Since $r^2|q|^2=|z|^2$, the first term is part of $L(z)$. By \eqref{eq:leibniz} applied to $q=\rho z$ and Proposition \ref{prop:weights} 
$$
r^2|D_h^nq|^2\le C\sum_{k=0}^{n}s^{2(n-k)}|D_h^{k}A_h^{n-k}z|^2,\qquad n=1,\dots,4.
$$
Multiplying by $s^{7-2n}$ and integrating, the bound $|A_h^kf|^2\le A_h^k(|f|^2)$ bounds the terms with $k\le3$ by $CL(z)$, and the term $k=4$, $n=4$ by $C\mathcal{F}(z)$. For the time derivative, \eqref{eq:Dt:product} and Lemma \ref{lem:time:conj} give $\trm r\,D_tq=D_tz+\big(\tau\varphi\,\trm(\theta')+\mathcal E\mathcal{O}(1)\big)\trp z$, and $|\tau\varphi\trm(\theta')|\le C\trm(s^2)$. Hence
\begin{align*}
\int_Q\trm(s^{-1}r^2)|D_tq|^2\le &2\mathcal{D}(z)+C\int_Q\trm(s^3)|\trp z|^2\le 2\mathcal{D}(z)+C\int_{Q_{\ast}}s^3|z|^2\\
&+C\Delta t\,s(T+\tfrac{\Delta t}2)^3\int_{\mathcal{M}}|\trp z(T)|^2,
\end{align*}
where we used \eqref{eq:theta:shift} to replace $\trm(s^3)$ by $\trp(s^3)$, and \eqref{time:integral}. The last term is bounded by $Ch^{-4}\mathcal{T}(z)$, since $(sh)^3\le8$ by \eqref{eq:s:bound} and $\Delta t\le1$.
\end{proof}

\subsection{End of the proof of Theorem \ref{thm:intro:carleman}}

We insert Propositions \ref{prop:spatial} and \ref{prop:time} into \eqref{eq:general estimate}, and use Lemma \ref{lem:Dt:absorb} and Lemma \ref{lem:estimate Rv}. This yields
$$
\|S_hz\|^2+\|P_hz\|^2+2c_0L(z)\le2\|\trm r\,\mathcal{P}^\star q\|^2+C\big(\varepsilon_0+\tau_0^{-1}\big)\big(L(z)+\mathcal{D}(z)\big)+C\mathcal{B}(z)+Ch^{-4}\mathcal{T}(z),
$$
the positive term $\displaystyle\frac{\Delta t}{2}\int|D_h^2D_tz|^2$ having absorbed the term $\frac{\Delta t}{8}\int|D_h^2D_tz|^2$ from Lemma \ref{lem:Dt:absorb}. By Lemmas \ref{lem:absorb} and \ref{lem:fourth:term}, for $\kappa\in(0,1]$, $\|S_hz\|^2+\|P_hz\|^2\ge\kappa\big(\frac18\mathcal{F}(z)+\frac12\mathcal{D}(z)\big)-C\kappa L(z)$. Choosing $\kappa$ small (depending on $c_0$), then $\varepsilon_0$ small and $\tau_0$ large, we obtain
$$
L(z)+\mathcal{D}(z)+\mathcal{F}(z)\le C\Big(\|\trm r\,\mathcal P^\star q\|^2_{L^2_h(Q)}+\mathcal{B}(z)+h^{-4}\mathcal{T}(z)\Big).
$$
Since $\displaystyle\mathcal T(z)=\int_{\mathcal{M}\times\partial\mathcal{N}}|\trp(rq)|^2$, Lemmas \ref{lem:boundary} and \ref{lem:primalvariable} yield \eqref{eq:intro:carleman}. \hfill$\square$

\section{\texorpdfstring{$\phi(h)$}{phi(h)}-null controllability}\label{sec:control}

\begin{lemma}\label{lem:energy}
Assume $\gamma\Delta t\le\frac14$. Then \eqref{sys:primal} and \eqref{sys:adjoint} are uniquely solvable, and every solution of \eqref{sys:adjoint} satisfies
\begin{equation}\label{ine:t:mas}
\|\trp q(0)\|^2_{L^2_h(\mathcal{M})}\le e^{4\gamma t}\,\|\trp q(t)\|^2_{L^2_h(\mathcal{M})}\qquad\text{for all }t\in\overline{\mathcal N} .
\end{equation}
\end{lemma}
\begin{proof}
Each implicit step of \eqref{sys:adjoint} amounts to inverting $I+\Delta t(D_h^4+a)$ with homogeneous boundary conditions, which is positive definite on $L^2_h(\mathcal{M})$ by \eqref{eq:green} since $\gamma\Delta t<1$. For $t\in\mathcal{N}$, the first equation of \eqref{sys:adjoint} reads $\trm q=\trp q-\Delta t(D_h^4\trm q+a\,\trm q)$; taking the inner product with $\trm q$ and using \eqref{eq:green},
$$
\|\trm q\|^2+\Delta t\int_{\overline{\mathcal M}}|D_h^2\trm q|^2+\Delta t\int_{\mathcal{M}}a|\trm q|^2=\langle\trp q,\trm q\rangle\le\|\trp q\|\,\|\trm q\| ,
$$
so $(1-\gamma\Delta t)\|\trm q(t)\|\le\|\trp q(t)\|$. Iterating, and using $(1-x)^{-1}\le e^{2x}$ for $0\le x\le\frac12$, we obtain \eqref{ine:t:mas}. 
\end{proof}

\begin{proof}[Proof of Theorem \ref{thm:intro:observability}]
We first choose the parameters. Let $\tau_0,\varepsilon_0,h_0,C$ be given by Theorem \ref{thm:intro:carleman}, and set
$$
\tau:=\tau_1T^2\Lambda,\qquad \delta:=\frac{\tau_1}{\varepsilon_0}\,\Lambda\,h^{m_\mu},\qquad \tau_1:=\max\{\tau_0,(2C)^{1/7}\},
$$
with $\Lambda=1+\frac1T+\gamma^{2/7}$. Then $\tau\ge\tau_0(T+T^2)$ and $\tau\ge(2C)^{1/7}T^2\gamma^{2/7}$. Moreover $\delta\le\frac12$ if $h\le h_1:=C_0\Lambda^{-\max\{1,5/\mu\}}$ with $C_0:=(\varepsilon_0/(2\tau_1))^{\max\{1,5/\mu\}}$, and
$$
\frac{\tau h}{\delta T^2}=\varepsilon_0h^{1-m_\mu}\le\varepsilon_0,\qquad
\frac{\tau^5\Delta t}{\delta^5T^{10}}=\varepsilon_0^5\frac{\Delta t}{h^{5m_\mu}}\le\varepsilon_0^5h^{\mu-5m_\mu}\le\varepsilon_0,\qquad
\frac{\tau}{\delta T^2}=\frac{\varepsilon_0}{h^{m_\mu}} .
$$
Hence Theorem \ref{thm:intro:carleman} applies.

We now apply the Carleman estimate. Since $\mathcal P^\star q=-a\,\trm q$, \eqref{eq:slicing} gives $\displaystyle\int_Q(\trm r)^2|\mathcal P^\star q|^2\le\gamma^2\int_{Q_{\ast}}r^2|q|^2$, and $s^7\ge(\tau/T^2)^7\ge2C\gamma^2$ on $[0,T]$. Thus this term is absorbed and
\begin{equation}\label{ine:obs}
\int_{Q_{\ast}}s^7r^2|q|^2\le C\bigg(\int_{\mathcal{N}^{\ast}}s^3r^2|D_h^2q|^2\Big|_{x=0}+\int_{\mathcal{N}^{\ast}}s\,r^2\big|_{x=0}|D_h^3q|^2\Big|_{x=h/2}+h^{-4}\int_{\mathcal{M}\times\partial\mathcal{N}}|\trp(rq)|^2\bigg).
\end{equation}

Next, we bound the left-hand side of \eqref{ine:obs} from below. For $t\in\mathcal{N}^{\ast}\cap(\frac T4,\frac{3T}4)$ we have $T^{-2}\le\theta\le\frac{16}{3T^2}$, so $s^7\ge1$ and $r^2=e^{-2s\varphi}\ge e^{-c_K\tau/T^2}$ with $c_K:=\frac{32}{3}K_0$. Every $t\in\mathcal{N}^{\ast}$ satisfies $q(t)=\trp q(t-\frac{\Delta t}2)$ with $t-\frac{\Delta t}{2}\in\overline{\mathcal N}$; hence, by \eqref{ine:t:mas} and since $\Delta t\le T/4$ implies that
$$
\int_{Q_{\ast}}s^7r^2|q|^2\ge\frac T4\,e^{-4\gamma T-c_K\tau/T^2}\,\|\trp q(0)\|^2_{L^2_h(\mathcal{M})} .
$$

We now bound the terms on the right-hand side of \eqref{ine:obs} from above. Since $\varphi\ge k_0$, $s^3r^2+s\,r^2\le\sup_{\sigma>0}(\sigma^3+\sigma)e^{-2k_0\sigma}=:C_k$. At the time boundary, the weight is evaluated at $t=\frac{\Delta t}2$ and $t=T+\frac{\Delta t}2$; the above choice of the parameters gives \eqref{eq:dt:small}, that is $\Delta t\le\delta T/16$, hence $\theta\ge(2\delta T^2)^{-1}$ at both points, so $r^2\le e^{-k_0\tau/(\delta T^2)}=e^{-k_0\varepsilon_0/h^{m_\mu}}$. Together with \eqref{ine:t:mas} for $t=T$,
$$
h^{-4}\int_{\mathcal{M}\times\partial\mathcal{N}}|\trp(rq)|^2\le2h^{-4}e^{4\gamma T}e^{-k_0\varepsilon_0/h^{m_\mu}}\|q_T\|^2_{L^2_h(\mathcal{M})} .
$$

Finally, combining the three previous estimates, with $\tau/T^2=\tau_1\Lambda$, we obtain
$$
\|\trp q(0)\|^2\le\frac{4C}{T}e^{8\gamma T+c_K\tau_1\Lambda}\bigg(C_k\int_{\mathcal{N}^{\ast}}\Big(|D_h^2q|^2\Big|_{x=0}+|D_h^3q|^2\Big|_{x=h/2}\Big)+2h^{-4}e^{-k_0\varepsilon_0/h^{m_\mu}}\|q_T\|^2\bigg).
$$
Since $T^{-1}\le e^{\Lambda}$, the prefactor is bounded by $C_{\rm obs}=e^{C(\Lambda+T\gamma)}$ (with $C_k$ absorbed); and, decreasing $h_0$ if necessary (depending on $\mu$), $2h^{-4}e^{-k_0\varepsilon_0/h^{m_\mu}}\le e^{-C_1/h^{m_\mu}}$ with $C_1:=k_0\varepsilon_0/2$. This proves \eqref{eq:intro:obs}.
\end{proof}

\begin{proof}[Proof of Theorem \ref{thm:intro:control}]
For $q_T\in L^2_h(\mathcal{M})$, let $q$ be the solution of \eqref{sys:adjoint}, $\displaystyle W(q):=\big(D_h^2q|_{x=0},\,D_h^3q|_{x=h/2}\big)$ on $\mathcal{N}^{\ast}$, $\|W\|^2:=\int_{\mathcal{N}^{\ast}}\big(|D_h^2q|_{x=0}|^2+|D_h^3q|_{x=h/2}|^2\big)$, and consider the penalized functional
$$
J(q_T):=\frac12\|W(q)\|^2+\frac{\phi(h)}2\|q_T\|^2_{L^2_h(\mathcal{M})}+\langle y_0,\trp q(0)\rangle .
$$
By \eqref{eq:intro:obs} and Young's inequality,
$$
|\langle y_0,\trp q(0)\rangle|\le\|y_0\|\big(C_{\rm obs}\|W\|^2+C_{\rm obs}\phi(h)\|q_T\|^2\big)^{1/2}\le C_{\rm obs}\|y_0\|^2+\frac14\|W\|^2+\frac{\phi(h)}4\|q_T\|^2,
$$
so that $J(q_T)\ge\frac14\|W\|^2+\frac{\phi(h)}4\|q_T\|^2-C_{\rm obs}\|y_0\|^2$. Since $L^2_h(\mathcal{M})$ is finite-dimensional and $J$ is continuous, strictly convex and coercive, it has a unique minimizer $q_T^{\rm opt}$, characterized by
\begin{equation}\label{eq:wopt}
\langle W^{\rm opt},W\rangle+\phi(h)\langle q_T^{\rm opt},q_T\rangle+\langle y_0,\trp q(0)\rangle=0\qquad\forall q_T\in L^2_h(\mathcal{M}),
\end{equation}
where $W^{\rm opt}:=W(q^{\rm opt})$. We define the controls on $\mathcal{N}^{\ast}$ by
$$
u_2:=D_h^2q^{\rm opt}\big|_{x=0},\qquad u_1:=-D_h^3q^{\rm opt}\big|_{x=h/2} .
$$
By Lemma \ref{lem:duality}, for every $q_T$, $$\langle y(T),q_T\rangle=\langle y_0,\trp q(0)\rangle+\langle W^{\rm opt},W\rangle=-\phi(h)\langle q_T^{\rm opt},q_T\rangle,$$ that is, $y(T)=-\phi(h)q_T^{\rm opt}$. Taking $q_T=q_T^{\rm opt}$ in \eqref{eq:wopt} and using the bound above,
$$
\|W^{\rm opt}\|^2+\phi(h)\|q_T^{\rm opt}\|^2=-\langle y_0,\trp q^{\rm opt}(0)\rangle\le C_{\rm obs}\|y_0\|^2+\frac14\|W^{\rm opt}\|^2+\frac{\phi(h)}4\|q_T^{\rm opt}\|^2,
$$
hence $\|u_1\|^2+\|u_2\|^2=\|W^{\rm opt}\|^2\le C_{\rm obs}\|y_0\|^2$ and
$$
\|y(T)\|=\sqrt{\phi(h)}\,\big(\phi(h)\|q_T^{\rm opt}\|^2\big)^{1/2}\le\sqrt{C_{\rm obs}}\,\sqrt{\phi(h)}\,\|y_0\| .
$$
\end{proof}

\section{Comments and conclusion}\label{sec:conclusion}

We have proved a Carleman estimate for the implicit Euler and centered finite-difference approximation of a one-dimensional fourth-order parabolic operator with clamped boundary conditions and boundary observation (Theorem \ref{thm:intro:carleman}). From it we deduced a relaxed observability inequality (Theorem \ref{thm:intro:observability}) and the $\phi$-null controllability of the fully discrete system with two boundary controls acting at $x=0$ (Theorem \ref{thm:intro:control}). The controls are bounded uniformly with respect to $h$ and $\Delta t$, and the final state is of order $\sqrt{\phi(h)}$ with $\phi(h)$ exponentially small. 

\par The present work belongs to the line of research  \cite{BHL10a,boyer-2014,CLTP:2022,GC-HS-2021,LMPZ:2023,AAP-2024}, in which the Carleman estimate is proved directly for the discrete operator, by conjugation with the weight and discrete integrations by parts. The closest fully discrete work outside the second-order case is \cite{kumar:2026}, which treats the implicit (backward Euler) finite-difference approximation of $\partial_t+\partial_x^3$ on $(0,L)$, with one Dirichlet control at $x=0$. As here, the weight is $e^{-s(t)\beta(x)}$ with $\theta$ regularized by a parameter $\delta>0$ (which allows the weight to be evaluated outside $[0,T]$), the observation is a discrete boundary derivative ($D_+D_-\varphi_0$ there, $D_h^2q(0)$ and $D_h^3q(\frac h2)$ here), and the estimate contains an additional negative power of the mesh size in front of the values at the initial and final times ($(\Delta x)^{-5/2}$ there, $h^{-4}$ here). The conditions of \cite[Theorem 1.6]{kumar:2026} are $\lambda^2\Delta x(\delta^2T^4)^{-1}\le\epsilon_0$ and $\lambda^3\Delta t(\delta^4T^6)^{-1}\le(\Delta x)^{\alpha}$ with $\alpha\ge1$; in terms of $S=\lambda(\delta T^2)^{-1}$, the first one reads $S^2\Delta x\le\epsilon_0$, which leads to the choice $\delta\approx(\Delta x)^{1/2}$ and to $\phi(\Delta x)=e^{-C_0/(\Delta x)^{1/2}}$ under $\Delta t<\delta(\Delta x)^{3/2+\alpha}$ \cite[Corollary 1.9, Theorem 1.11]{kumar:2026}. In the present fourth-order setting the spatial condition is linear in $S$, namely $Sh\le\varepsilon_0$, which is what allows $\delta\approx h^{m_\mu}$ and the rate $e^{-C_1/h}$ for $\mu\ge5$.
\par Several questions remain open. Extending the transfer principle of \cite{LW:2026} to fourth-order operators with boundary observation, and comparing the resulting range with \eqref{eq:intro:conditions}, is a natural next step. Removing the coupling $\Delta t\le h^{\mu}$ between the discretization parameters, proving $\phi$-null controllability with only one of the two boundary controls, and extending Theorem \ref{thm:intro:carleman} to several space dimensions with the estimates of \cite{AAP-2024} are further directions. In \cite{guerrero-2019}, Guerrero and Kassab proved a global Carleman estimate for the fourth-order parabolic operator in a bounded domain $\Omega\subset\mathbb R^N$, $N\ge2$, with an observation localized in an open subset $\omega\subset\Omega$, and deduced the null controllability of the corresponding equation. Since the observation is distributed in $\omega$, this is the setting of \cite{LW:2026}, and it is natural to ask whether the transfer principle of \cite{LW:2026}, with the estimate of \cite{guerrero-2019} playing the role of the continuous Carleman inequality (A4), yields a fully discrete Carleman estimate and the $\phi$-null controllability of a fully discrete fourth-order equation in dimension $N\ge2$. 
Such a result would complement Theorem \ref{thm:intro:control}, which is one-dimensional and uses boundary controls. A direct proof along the lines of the present paper would require the multidimensional weight calculus of \cite{AAP-2024} for operators of order four, together with a fully discrete version of the computations of \cite{guerrero-2019}.

\appendix
\section{Cross terms involving the discrete time derivative}\label{estimaciones:carleman}

In this appendix we estimate $I_{11}+I_{41}$, $I_{21}$ and $I_{31}$, under the assumptions of Section \ref{sec:proof:carleman}. Recall that $w=D_tz$ vanishes on $(\partial\mathcal{M}\cup\partial\overline{\mathcal M})\times\mathcal{N}$, and that, by \eqref{eq:s:bound}, $s\le2S$ and $sh\le2\varepsilon_0$ at every time where the weights are evaluated. Time-boundary terms are bounded by an inverse inequality. For $u$ vanishing on $\partial\mathcal{M}\cup\partial\overline{\mathcal M}$ and $x\in\overline{\mathcal M}$, the Cauchy--Schwarz inequality  gives
$$
h^4|D_h^2u(x)|^2=|u(x+h)-2u(x)+u(x-h)|^2\le6\big(|u(x+h)|^2+|u(x)|^2+|u(x-h)|^2\big),
$$
and, since each value of $u$ appears at most three times (twice for $h^2|D_hu|^2$),
\begin{equation}\label{ine:useful}
\int_{\overline{\mathcal M}}h^4|D_h^2u|^2\le18\int_{\mathcal{M}}|u|^2,\qquad \int_{\widehat{\mathcal M}^{\ast}}h^2|D_hu|^2\le4\int_{\mathcal{M}}|u|^2 .
\end{equation}
Together with $|A_h^kz|^2\le A_h^k(|z|^2)$ and $s\le2S$, this gives for instance
\begin{equation}\label{eq:tb}
\int_{\overline{\mathcal M}}s^4|A_h^2\trp z|^2+\int_{\mathcal M^{\ast}}s^2|D_h\trp z|^2+\int_{\overline{\mathcal M}}|D_h^2\trp z|^2\le C\big((sh)^4+(sh)^2+1\big)h^{-4}\int_{\mathcal{M}}|\trp z|^2 .
\end{equation}

\subsection{Estimate of \texorpdfstring{$I_{11}$ and $I_{41}$}{I11 and I41}}
\begin{lemma}\label{lem:estimate:I11}
There exists $C_1>0$ such that
$$
I_{11}+I_{41}\ge-C\int_{Q^{\ast}}T\theta s^2|D_hz|^2-C_1\Delta t\int_{\mathcal M^{\ast}\times\mathcal{N}}\trm(s^2)|D_hD_tz|^2-C\varepsilon_0^2\big(L(z)+\mathcal{D}(z)\big)-C\varepsilon_0^2h^{-4}\mathcal T(z).
$$
\end{lemma}
\begin{proof}
Set $\alpha_{11}:=rA_h^2D_h^2\rho$, so that $\displaystyle          I_{11}=-6\int_Q \trm\alpha_{11}\,D_h^2A_h^2\trm z\, D_{t}z$. By the first identity in \eqref{eq:ipp}, with $D_{t}z=0$ on $\partial\mathcal{M}$, we have
$$I_{11}=6\int_{\mathcal{M}^{\ast}\times \mathcal{N}} D_{h}A_{h}^{2}\trm z\,D_{h}(\trm\alpha_{11} D_{t}z).$$
Then by the second identity in \eqref{eq:ipp}, with $D_{t}z=0$ on $\partial\overline{\mathcal{M}}$, it follows that
$$
I_{11}=6\int_{\overline{\mathcal M}\times\mathcal{N}}D_hA_h\trm z\;D_hA_h(\trm\alpha_{11} D_{t}z).
$$
By Lemma \ref{pro:derivadapromedio}, $$D_hA_h(\trm\alpha_{11} D_{t}z)=D_hA_h\trm\alpha_{11}\,(D_{t}z+\frac{h^2}2D_h^2D_{t}z)+(A_h^2(\trm\alpha_{11})+\frac{h^2}4D_h^2(\trm\alpha_{11}))\,D_hA_hD_{t}z.$$ Since $S_4z=6\trm(D_hA_h\alpha_{11}\,A_hD_hz)$, we have $\displaystyle I_{41}=-6\int_QD_hA_h\trm\alpha_{11}\,A_hD_h\trm z\,D_{t}z$, and the integrand vanishes on $\partial\mathcal{M}$; hence the first term cancels and
$$
I_{11}+I_{41}=3h^2\int_{\overline{\mathcal M}\times\mathcal{N}}D_hA_h\trm\alpha_{11}\,A_hD_h\trm z\;D_h^2D_{t}z+6\int_{\overline{\mathcal M}\times\mathcal{N}}\trm\beta\;D_hA_hD_{t}z\,A_hD_h\trm z=:K_1+K_2,
$$
with $\beta:=A_h^2\alpha_{11}+\frac{h^2}4D_h^2\alpha_{11}$. By Proposition \ref{prop:weights}, $\beta=r\partial_x^2\rho+s^2\mathcal{O}((sh)^2)=s^2(\partial_x\varphi)^2+s\mathcal{O}(1)+s^2\mathcal{O}((sh)^2)$, so $0\le\beta\le Cs^2$.

\par Let us focus on the term $K_2$. By \eqref{eq:identity:square}, $D_hA_hD_{t}z A_hD_h\trm z=\frac12D_t|A_hD_hz|^2-\frac{\Delta t}2|D_tA_hD_hz|^2$. Using time-discrete integration by parts \eqref{eq:tau:mas},
$$
K_2=-3\int_{\overline{\mathcal M}\times\mathcal{N}^{\ast}}D_t(\trm\beta)\,|A_hD_hz|^2+3\int_{\overline{\mathcal M}\times\partial\mathcal{N}}\trm\beta\;\trp|A_hD_hz|^2\,n-3\Delta t\int_{\overline{\mathcal M}\times\mathcal{N}}\trm\beta\,|D_tA_hD_hz|^2 .
$$
By Lemma \ref{lem:Dt:outer}, $D_t(\trm\beta)=T\theta s^2\mathcal{O}(1)$; with $|A_hD_hz|^2\le A_h(|D_hz|^2)$ and \eqref{eq:ipp} with $\mathcal W=\overline{\mathcal M}$, the first term is bounded below by $\displaystyle-C\int_{Q^{\ast}}T\theta s^2|D_hz|^2$. The boundary term at $t=T$ is nonnegative; the one at $t=0$ is bounded below by $-C\varepsilon_0^2h^{-4}\mathcal{T}(z)$ after using the same strategy to obtain \eqref{eq:tb}. In the last term, $|D_{t}A_{h}D_{h}z|^2\le A_{h}(|D_{h}D_{t}z|^2)$ and $D_{h}D_{t}z=0$ on $\partial\widehat{\mathcal M}^{\ast}$, so it is bounded below by $\displaystyle-C_1\Delta t\int_{\mathcal M^{\ast}\times\mathcal{N}}\trm(s^2)|D_hD_{t}z|^2$.

Finally, for the term $K_1$, it follows by \eqref{eq:ipp} with $\mathcal W=\overline{\mathcal M}$ and with $\mathcal W=\mathcal{M}$, using $D_{h}D_{t}z=0$ on $\partial\widehat{\mathcal M}^{\ast}$ and $D_{t}z=0$ on $\partial\mathcal{M}$, $$\int_{\overline{\mathcal M}}D_{h}A_{h}\trm\alpha_{11}\,A_{h}D_{h}\trm z\,D_{h}^{2}D_{t}z=\int_{\mathcal{M}}D_{t}z\,D_{h}^{2}(D_{h}A_{h}\trm\alpha_{11}\,A_{h}D_{h}\trm z)$$
By \eqref{eq:leibniz} and Proposition \ref{prop:weights} we have $D_{h}^{k}A_{h}^{l} \trm\alpha_{11}=s^2\mathcal{O}(1)$,
$$|D_h^2(D_{h}A_{h}\trm\alpha_{11}A_{h}D_{h}\trm z)|\le Cs^2\big(|A_h^3D_h\trm z|+|A_h^2D_h^2\trm z|+|A_hD_h^3\trm z|\big).$$ Hence, by Young's inequality, $|A_{h}^{k}\trm z|^2\le A_h^k(|\trm z|^2)$ and \eqref{eq:slicing},
$$
|K_1|\le C(sh)^2\Big(\mathcal{D}(z)+\int_{Q^{\ast}}s|D_hz|^2+\int_{\overline{Q}_{\ast}}s|D_h^2z|^2+\int_{Q^{\ast}}s|D_h^3z|^2\Big)\le C\varepsilon_0^2\big(\mathcal{D}(z)+L(z)\big).
$$
\end{proof}
\subsection{Estimate of \texorpdfstring{$I_{21}$}{I21}}
\begin{lemma}\label{lem:estimate:I21}
We have
$$
I_{21}\ge-C\int_{Q_{\ast}}T\theta s^4|z|^2-C\varepsilon_0\big(L(z)+\mathcal{D}(z)\big)-C\varepsilon_0^4h^{-4}\mathcal{T}(z).
$$
\end{lemma}
\begin{proof}
Set $\alpha_{21}:=rD_h^4\rho$, so that $\displaystyle I_{21}=-\int_Q \trm\alpha_{21}\,A_h^4\trm z\;D_{t}z$. Using the averaging formulas of Lemma \ref{lem:ipp} twice, with $\trm\alpha_{21}\,D_{t}z=0$ on $\partial\mathcal{M}$ and $A_h(\trm\alpha_{21}\,D_{t}z)=0$ on $\partial\widehat{\mathcal M}^{\ast}$,
$$
I_{21}=-\int_{\overline{\mathcal M}\times\mathcal{N}}A_h^2(\trm\alpha_{21}\,D_{t}z)\,A_h^2\trm z .
$$
By Lemma \ref{pro:derivadapromedio}, $$A_h^2(\trm\alpha_{21}D_{t}z)=A_h^2(\trm\alpha_{21})\,A_h^2(D_{t}z)+\frac{h^2}2D_hA_h(\trm\alpha_{21})\,D_hA_h(D_{t}z)+\frac{h^4}{16}D_h^2(\trm\alpha_{21})\,D_h^2D_{t}z,$$ and we split $I_{21}=I_{21}^{(a)}+I_{21}^{(b)}+I_{21}^{(c)}$ accordingly.

For the term $I_{21}^{(a)}$, we use \eqref{eq:identity:square} to write $A_h^2(D_{t}z)\,A_h^2\trm z=\frac{1}{2}D_t|A_h^2z|^2-\frac{\Delta t}2|D_tA_h^2z|^2$, and $A_h^2(\trm\alpha_{21})=\trm(r\partial_x^4\rho)+s^4\mathcal{O}((sh)^2)\ge0$ by Proposition \ref{prop:weights}. We drop the nonnegative term $\displaystyle\frac{\Delta t}2\int_{\overline{\mathcal{M}}\times\mathcal{N}} A_h^2(\trm\alpha_{21})\,|D_tA_h^2z|^2$ and apply time-discrete integration by parts \eqref{eq:tau:mas}, which gives
$$
I_{21}^{(a)}\ge\frac12\int_{\overline{\mathcal M}\times\mathcal{N}^{\ast}}D_t(A_h^2(\trm\alpha_{21}))\,|A_h^2z|^2-\frac12\int_{\overline{\mathcal M}\times\partial\mathcal{N}}A_h^2(\trm\alpha_{21})\;\trp|A_{h}^{2}z|^2\,n .
$$
By Lemma \ref{lem:Dt:outer}, $D_t(A_h^2(\trm\alpha_{21}))=T\theta s^4\mathcal{O}(1)$, and $\displaystyle\int_{\overline{\mathcal M}}|A_h^2z|^2\le\int_{\mathcal{M}}|z|^2$; the second integral in the right-hand side above, which are time boundary terms, are bounded by $C\varepsilon_0^4h^{-4}\mathcal{T}(z)$.

Now, let us focus on the terms $I_{21}^{(b)}$ and $I_{21}^{(c)}$. By Proposition \ref{prop:weights}, $D_{h}A_{h}\trm\alpha_{21}$ and $D_h^2\trm\alpha_{21}$ are $s^4\mathcal{O}(1)$. Since $h|D_{h}A_{h}D_{t}z(x)|\le|A_{h}D_{t}z(x+\frac h2)|+|A_{h}D_{t}z(x-\frac h2)|\le |D_{t}z(x+h)|+2|D_{t}z(x)|+|D_{t}z(x-h)|$ and $h^2|D_{h}^{2}D_{t}(x)|\le|D_{t}(x+h)|+2|D_{t}(x)|+|D_{t}(x-h)|$,
$$
|I_{21}^{(b)}|+|I_{21}^{(c)}|\le C\int_{\overline{\mathcal M}\times\mathcal{N}}\big(h\,\trm(s^4)+h^2\trm(s^4)\big)\,|D_{t}z|\;|A_h^2\trm z|.
$$
 Writing $hs^4\,D_{t}z\,A_{h}^{2}\trm z=(sh)\,(s^{-1/2}D_{t}z)(s^{7/2}A_{h}^{2}\trm z)$ and using Young's inequality, we have
 $$
|I_{21}^{(b)}|+|I_{21}^{(c)}|\le C\int_{\overline{\mathcal M}\times\mathcal{N}}h\,\trm(s^{-1})|D_{t}z|^{2}+\trm(sh)\trm(s^{7})\,|A_h^2\trm z|^{2}.
$$
Using $|A_h^2\trm z|^2\le A_h^2(|\trm z|^2)$ and \eqref{eq:slicing}, we obtain $\displaystyle|I_{21}^{(b)}|+|I_{21}^{(c)}|\le C\varepsilon_0\left(\mathcal{D}(z)+\int_{Q_{\ast}}s^7|z|^2\right)$.
\end{proof}

\subsection{Estimate of \texorpdfstring{$I_{31}$}{I31}}

\begin{lemma}\label{lem:estimate:I31}
We have
$$
I_{31}\ge\frac{\Delta t}4\int_{\overline{\mathcal M}\times\mathcal{N}}|D_h^2D_tz|^2-C\int_{\overline{Q}_{\ast}}T\theta|D_h^2z|^2-C\varepsilon_0^2\left(\mathcal{D}(z)+\int_{\overline{Q}_{\ast}}s|D_h^2z|^2+\int_{Q^{\ast}}s|D_h^3z|^2\right)-18h^{-4}\mathcal T(z).
$$
\end{lemma}
\begin{proof}
Set $\alpha_{31}:=rA_{h}^{4}\rho$, so that $\displaystyle I_{31}=-\int_{Q} \trm\alpha_{31}\,D_{t}z\,D_{h}^{4}\trm z$. Notice that the functions $\trm\alpha_{31}D_{t}z$ and $\trm z$ vanish on the boundary $\partial\mathcal{M}$ and also their spatial derivative on $\partial\widehat{\mathcal M}^{\ast}$. Thus, applying \eqref{eq:green}; then, using \eqref{eq:leibniz} with $n=2$, we obtain
$$
I_{31}=-\int_{\overline{\mathcal M}\times\mathcal{N}}D_h^2(\trm \alpha_{31}D_{t}z)\,D_h^2\trm z=I_{31}^{(a)}+I_{31}^{(b)}+I_{31}^{(c)},
$$
with
$$
\begin{gathered}
I_{31}^{(a)}:=-\int_{\overline{\mathcal M}\times\mathcal{N}}A_{h}^{2}\trm\alpha_{31}\,D_h^2D_{t}z\,D_h^2\trm z,\qquad I_{31}^{(b)}:=-2\int_{\overline{\mathcal M}\times\mathcal{N}}D_hA_h\trm\alpha_{31}\,A_hD_hD_{t}z\,D_h^2\trm z,\\
I_{31}^{(c)}:=-\int_{\overline{\mathcal M}\times\mathcal{N}}D_h^2\trm\alpha_{31}\,A_h^2D_{t}z\,D_h^2\trm z .
\end{gathered}
$$
By Proposition \ref{prop:weights}, $A_h^2\trm\alpha_{31}=1+\mathcal{O}((sh)^2)\in[\frac12,2]$, while $D_h^kA_h^l\trm\alpha_{31}=\mathcal{O}((sh)^2)$ for $k\ge1$; by Lemma \ref{lem:Dt:outer}, $D_t(A_{h}^{2}\trm\alpha_{31})=T\theta\,\mathcal{O}(1)$.

 Since $D_h^2w\,D_h^2\trm z=\frac12D_t|D_h^2z|^2-\frac{\Delta t}2|D_tD_h^2z|^2$, for the term $I_{31}^{(a)}$, \eqref{eq:tau:mas} gives
$$
I_{31}^{(a)}=\frac12\int_{\overline{\mathcal M}\times\mathcal{N}^{\ast}}D_{t}(A_h^2\trm\alpha_{31})\,|D_{h}^{2}z|^2-\frac12\int_{\overline{\mathcal M}\times\partial\mathcal{N}}A_{h}^{2}\trm\alpha_{31}\,\trp|D_{h}^{2}z|^2\,n+\frac{\Delta t}2\int_{\overline{\mathcal M}\times\mathcal{N}}A_{h}^{2}\trm\alpha_{31}\,|D_{h}^{2}D_{t}z|^2.
$$
The first term is bounded below by $\displaystyle-C\int_{\overline{Q}_{\ast}}T\theta|D_h^2z|^2$, the last one by $\displaystyle\frac{\Delta t}4\int_{\overline{\mathcal{M}}\times\mathcal{N}}|D_h^2D_tz|^2$. \break The boundary term at $t=0$ is nonnegative, and the one at $t=T$ is bounded below by \break $\displaystyle-\int_{\overline{\mathcal M}}|D_h^2\trp z(T)|^2\ge-18h^{-4}\int_{\mathcal{M}}|\trp z(T)|^2$ by \eqref{ine:useful}.

Let us focus on the term $I_{31}^{(b)}$. Applying the average discrete integration by parts, the second identity in \eqref{eq:ipp}, and noticing that $D_{h}D_{t}z$ has vanishing boundary term, and the nodes of $\partial\widehat{\mathcal M}^{\ast}$ do not contribute to the integral over $\widehat{\mathcal M}^{\ast}$; it follows that
$$
I_{31}^{(b)}=-2\int_{\mathcal M^{\ast}}D_hD_{t}z\,A_h(D_{h}A_{h}\trm\alpha_{31}D_{h}^{2}\trm z).
$$
Since $D_{t}z=0$ on $\partial\mathcal{M}$, the first identity in \eqref{eq:ipp} with $\mathcal W=\mathcal{M}$, gives
\begin{equation}\label{eq:I31b:ipp}
I_{31}^{(b)}=2\int_{\mathcal{M}\times\mathcal{N}}D_{t}z\,D_hA_h(D_{h}A_{h}\trm \alpha_{31}D_{h}^{2}\trm z).
\end{equation}
We now expand $D_{h}A_{h}(D_{h}A_{h}\trm \alpha_{31}D_{h}^{2}\trm z)$. Applying Lemma \ref{pro:derivadapromedio} twice, first the second identity in \eqref{eq:difference:product} and then the first one to each of the two resulting products, we obtain, in $\mathcal{M}$,
\begin{align*}
D_hA_h(D_{h}A_{h}\trm\alpha_{31}D_{h}^{2}\trm z)=&D_{h}^{2}A_{h}^{2}\trm \alpha_{31}\,A_{h}^{2}D_{h}^{2}\trm z+A_{h}^{3}D_{h}\trm\alpha_{31}\,A_{h}D_{h}^{3}\trm z\\
&+\frac{h^2}{4}\big(D_{h}^{3}A_{h}\trm\alpha_{31}\,A_{h}D_{h}^{3}\trm z+D_{h}^{2}A_{h}^{2}\trm\alpha_{31}\,D_{h}^{4}\trm z\big).
\end{align*}
Proposition \ref{prop:weights} gives $\big(D_h^{k}A_h^{l}\trm\alpha_{31}\big)=\mathcal{O}((sh)^2)$ for all $k,l\in\mathbb{N}$. Moreover, for $x\in\mathcal{M}$,
$$
h^2|D_h^4\trm z(x)|\le|D_h^2\trm z(x+h)|+2|D_h^2\trm z(x)|+|D_h^2\trm z(x-h)|.
$$
Using also $h\le1$, we conclude that, in $\mathcal{M}$,
\begin{equation}\label{eq:I31b:bound}
|D_hA_h(D_{h}A_{h}\trm\alpha_{31}D_{h}^{2}\trm z)|\le C(sh)^2\Big(|A_h^2D_h^2\trm z|+|A_hD_h^3\trm z|+\sum_{j=-1}^{1}|D_h^2\trm z(\cdot+jh)|\Big).
\end{equation}
We now insert \eqref{eq:I31b:bound} in \eqref{eq:I31b:ipp} and apply Young's inequality, $(sh)^2\le4\varepsilon_0^2$, and therefore
$$
|I_{31}^{(b)}|\le C\varepsilon_0^2\int_{\mathcal{M}\times\mathcal{N}}\trm(s^{-1})|D_{t}z|^2+C\varepsilon_0^2\int_{\mathcal{M}\times\mathcal{N}}\trm(s)\Big(|A_h^2D_h^2\trm z|^2+|A_hD_h^3\trm z|^2+\sum_{j=-1}^{1}|D_h^2\trm z(\cdot+jh)|^2\Big).
$$
The first integral is $\mathcal{D}(z)$. For the second one, the weight depends only on time, and $|A_h^kf|^2\le A_h^k(|f|^2)$. Moreover, for $g\ge0$, writing $A_h^2g(x)=\frac14\big(g(x-h)+2g(x)+g(x+h)\big)$ and $A_hg(x)=\frac12\big(g(x-\frac h2)+g(x+\frac h2)\big)$ and summing over $x\in\mathcal{M}$, we have $$ \int_{\mathcal{M}}A_h^2g\le\int_{\overline{\mathcal M}}g, \quad\int_{\mathcal{M}}\sum_{j=-1}^1g(\cdot+jh)\le3\int_{\overline{\mathcal M}}g$$ and $\displaystyle\int_{\mathcal{M}}A_hg\le\int_{\mathcal M^{\ast}}g$. Hence
$$
\int_{\mathcal{M}}\Big(|A_h^2D_h^2\trm z|^2+|A_hD_h^3\trm z|^2+\sum_{j=-1}^{1}|D_h^2\trm z(\cdot+jh)|^2\Big)\le4\int_{\overline{\mathcal M}}|D_h^2\trm z|^2+\int_{\mathcal M^{\ast}}|D_h^3\trm z|^2 .
$$
Finally, \eqref{eq:slicing} turns the time integrals over $\mathcal{N}$ of functions of the form $\trm(\cdot)$ into integrals over $\mathcal{N}^{\ast}$, and we obtain
$$
|I_{31}^{(b)}|\le C\varepsilon_0^2\Big(D(z)+\int_{\overline{Q}_{\ast}}s|D_h^2z|^2+\int_{Q^{\ast}}s|D_h^3z|^2\Big).
$$

Since $D_{h}^{2}(\trm\alpha_{31})=\mathcal{O}((sh)^2)$, Young's inequality gives, at each $t\in\mathcal{N}$,
$$
|I_{31}^{(c)}|\le C(sh)^2\int_{\overline{\mathcal M}}\Big(\trm(s^{-1})|A_h^2w|^2+\trm(s)|D_h^2\trm z|^2\Big).
$$
As $D_{t}z=0$ on $\partial\mathcal{M}\cup\partial\overline{\mathcal M}$, we have $\displaystyle\int_{\overline{\mathcal M}}|A_h^2D_{t}z|^2\le\int_{\overline{\mathcal M}}A_h^2(|D_{t}z|^2)\le\int_{\mathcal{M}}|D_{t}z|^2$. Using $(sh)^2\le4\varepsilon_0^2$ and \eqref{eq:slicing} once more, we conclude that $\displaystyle|I_{31}^{(c)}|\le C\varepsilon_0^2\big(\mathcal{D}(z)+\int_{\overline{Q}_{\ast}}s|D_h^2z|^2\big)$.
\end{proof}

\bibliographystyle{acm}
\bibliography{references}

@article{sivashinsky-1977,
  AUTHOR  = {Sivashinsky, G. I.},
  TITLE   = {Nonlinear analysis of hydrodynamic instability in laminar flames. {I}. {D}erivation of basic equations},
  JOURNAL = {Acta Astronaut.},
  VOLUME  = {4},
  NUMBER  = {11-12},
  PAGES   = {1177--1206},
  YEAR    = {1977},
  DOI     = {10.1016/0094-5765(77)90096-0},
}

@article{kuramoto-tsuzuki-1975,
  AUTHOR  = {Kuramoto, Y. and Tsuzuki, T.},
  TITLE   = {On the formation of dissipative structures in reaction-diffusion systems: reductive perturbation approach},
  JOURNAL = {Progr. Theoret. Phys.},
  VOLUME  = {54},
  NUMBER  = {3},
  PAGES   = {687--699},
  YEAR    = {1975},
}

@article{kuramoto-tsuzuki-1976,
  AUTHOR  = {Kuramoto, Y. and Tsuzuki, T.},
  TITLE   = {Persistent propagation of concentration waves in dissipative media far from thermal equilibrium},
  JOURNAL = {Progr. Theoret. Phys.},
  VOLUME  = {55},
  NUMBER  = {2},
  PAGES   = {356--369},
  YEAR    = {1976},
  DOI     = {10.1143/PTP.55.356},
}

@article{takahashi-2017,
    AUTHOR = {Takahashi, Tak\'eo},
     TITLE = {Boundary local null-controllability of the
              {K}uramoto-{S}ivashinsky equation},
   JOURNAL = {Math. Control Signals Systems},
  FJOURNAL = {Mathematics of Control, Signals, and Systems},
    VOLUME = {29},
      YEAR = {2017},
    NUMBER = {1},
     PAGES = {Art. 2, 21},
      ISSN = {0932-4194,1435-568X},
   MRCLASS = {93B05 (35K35 35K59 93C10)},
  MRNUMBER = {3590382},
MRREVIEWER = {Nicol\'as\ Carre\~no},
       DOI = {10.1007/s00498-016-0182-5},
       URL = {https://doi-org.bibliotecavirtual.uis.edu.co/10.1007/s00498-016-0182-5},
}

@misc{HSM:2025,
  AUTHOR        = {Hern\'{a}ndez-Santamar\'{\i}a, V. and Majumdar, S.},
  TITLE         = {On the controllability of the {K}uramoto--{S}ivashinsky equation on multi-dimensional cylindrical domains},
  YEAR          = {2025},
  EPRINT        = {2508.00812},
  ARCHIVEPREFIX = {arXiv},
  NOTE          = {arXiv:2508.00812},
}

@article{LRR:2020,
    AUTHOR = {Le Rousseau, J\'er\^ome and Robbiano, Luc},
     TITLE = {Spectral inequality and resolvent estimate for the
              bi-{L}aplace operator},
   JOURNAL = {J. Eur. Math. Soc. (JEMS)},
  FJOURNAL = {Journal of the European Mathematical Society (JEMS)},
    VOLUME = {22},
      YEAR = {2020},
    NUMBER = {4},
     PAGES = {1003--1094},
      ISSN = {1435-9855,1435-9863},
   MRCLASS = {58J50 (35B45 35J30 35L35 35P15 35S15 93B05 93D15)},
  MRNUMBER = {4071321},
MRREVIEWER = {Marcos\ L\'opez-Garc\'ia},
       DOI = {10.4171/JEMS/939},
       URL = {},
}

@article{CGM:2017,
    AUTHOR = {Cerpa, Eduardo and Guzm\'an, Patricio and Mercado, Alberto},
     TITLE = {On the control of the linear {K}uramoto-{S}ivashinsky
              equation},
   JOURNAL = {ESAIM Control Optim. Calc. Var.},
  FJOURNAL = {ESAIM. Control, Optimisation and Calculus of Variations},
    VOLUME = {23},
      YEAR = {2017},
    NUMBER = {1},
     PAGES = {165--194},
      ISSN = {1292-8119,1262-3377},
   MRCLASS = {93B05 (35Q53 93C20)},
  MRNUMBER = {3601020},
MRREVIEWER = {Peng\ Gao},
       DOI = {10.1051/cocv/2015044},
       URL = {},
}

@misc{LW:2026,
  AUTHOR        = {L\"{u}, Qi and Wang, Yu},
  TITLE         = {From continuous to fully discrete {C}arleman estimates: a transfer principle for parabolic schemes},
  YEAR          = {2026},
  EPRINT        = {2610.02983},
  ARCHIVEPREFIX = {arXiv},
  NOTE          = {arXiv:2610.02983},
}

@article {kumar:2026,
    AUTHOR = {Kumar, Manish},
     TITLE = {Uniform boundary controllability of a fully discrete linear
              {K}d{V} equation {\it via} discrete {C}arleman estimates},
   JOURNAL = {ESAIM Control Optim. Calc. Var.},
  FJOURNAL = {ESAIM. Control, Optimisation and Calculus of Variations},
    VOLUME = {32},
      YEAR = {2026},
     PAGES = {Paper No. 69, 46},
      ISSN = {1292-8119,1262-3377},
   MRCLASS = {93B05 (35Q53 65M06)},
  MRNUMBER = {5128138},
       DOI = {10.1051/cocv/2026050},
       URL = {https://doi.org/10.1051/cocv/2026050},
}

@article {WZ:2024,
    AUTHOR = {Wang, Yu and Zhao, Qingmei},
     TITLE = {The {$\phi $}-{N}ull {C}ontrollability for {S}tochastic
              {F}ourth {O}rder {S}emi-{D}iscrete {P}arabolic {E}quations},
   JOURNAL = {Appl. Math. Optim.},
  FJOURNAL = {Applied Mathematics and Optimization},
    VOLUME = {95},
      YEAR = {2027},
    NUMBER = {1},
     PAGES = {Paper No. 3},
      ISSN = {0095-4616,1432-0606},
   MRCLASS = {93B05 (93B07 93C20)},
  MRNUMBER = {5149766},
       DOI = {10.1007/s00245-026-10527-8},
       URL = {https://doi.org/10.1007/s00245-026-10527-8},
}

@article {HS:2023,
    AUTHOR = {Hern\'andez-Santamar\'ia, V\'ictor},
     TITLE = {Controllability of a simplified time-discrete stabilized
              {K}uramoto-{S}ivashinsky system},
   JOURNAL = {Evol. Equ. Control Theory},
  FJOURNAL = {Evolution Equations and Control Theory},
    VOLUME = {12},
      YEAR = {2023},
    NUMBER = {2},
     PAGES = {459--501},
      ISSN = {2163-2472,2163-2480},
   MRCLASS = {93B05 (35K52 93C55)},
  MRNUMBER = {4520409},
MRREVIEWER = {V\~u{} Ngoc Ph\'at},
       DOI = {10.3934/eect.2022038},
       URL = {https://doi.org/10.3934/eect.2022038},
}

@article {zuazua:2005,
    AUTHOR = {Zuazua, Enrique},
     TITLE = {Propagation, observation, and control of waves approximated by
              finite difference methods},
   JOURNAL = {SIAM Rev.},
  FJOURNAL = {SIAM Review},
    VOLUME = {47},
      YEAR = {2005},
    NUMBER = {2},
     PAGES = {197--243},
      ISSN = {1095-7200,0036-1445},
   MRCLASS = {92C20 (35L05 65M06 65M30 93B05 93B07)},
  MRNUMBER = {2179896},
MRREVIEWER = {Fredi\ Tr\"oltzsch},
       DOI = {10.1137/S0036144503432862},
       URL = {https://doi.org/10.1137/S0036144503432862},
}

@article{AAP-2024,
  doi = {10.48550/ARXIV.2412.19892},
  
  url = {https://arxiv.org/abs/2412.19892},
  
  author = {Pérez, Ariel A.},
  JOURNAL = {arXiv:2412.19892 },
  
  title = {Asymptotic behavior of {C}arleman weight functions and application to controllability},
  
  publisher = {arXiv},
  
  year = {2026},
  
  copyright = {arXiv.org perpetual, non-exclusive license}
}

@article {LMPZ:2023,
    AUTHOR = {Lecaros, Rodrigo and Morales, Roberto and P\'erez, Ariel and
              Zamorano, Sebasti\'an},
     TITLE = {Discrete {C}arleman estimates and application to
              controllability for a fully-discrete parabolic operator with
              dynamic boundary conditions},
   JOURNAL = {J. Differential Equations},
  FJOURNAL = {Journal of Differential Equations},
    VOLUME = {365},
      YEAR = {2023},
     PAGES = {832--881},
      ISSN = {0022-0396,1090-2732},
   MRCLASS = {93B05 (65M06 93B07)},
  MRNUMBER = {4593117},
MRREVIEWER = {Salah-Eddine\ Chorfi},
       DOI = {10.1016/j.jde.2023.05.014},
       URL = {},
}

@article {LOP:2021,
    AUTHOR = {Lecaros, Rodrigo and Ortega, Jaime H. and P\'{e}rez, Ariel},
     TITLE = {Stability estimate for the semi-discrete linearized
              {B}enjamin-{B}ona-{M}ahony equation},
   JOURNAL = {ESAIM Control Optim. Calc. Var.},
  FJOURNAL = {ESAIM. Control, Optimisation and Calculus of Variations},
    VOLUME = {27},
      YEAR = {2021},
     PAGES = {Paper No. 93, 30},
      ISSN = {1292-8119},
   MRCLASS = {35B60 (35L05 35Q35 35R45 65M06)},
  MRNUMBER = {4313822},
       DOI = {10.1051/cocv/2021087},
       }

@article {GC-HS-2021,
    AUTHOR = {Gonz\'{a}lez Casanova, Pedro and Hern\'{a}ndez-Santamar\'{\i}a, V\'{\i}ctor},
     TITLE = {Carleman estimates and controllability results for fully
              discrete approximations of 1{D} parabolic equations},
   JOURNAL = {Adv. Comput. Math.},
  FJOURNAL = {Advances in Computational Mathematics},
    VOLUME = {47},
      YEAR = {2021},
    NUMBER = {5},
     PAGES = {Paper No. 72, 71},
      ISSN = {1019-7168},
   MRCLASS = {93B05 (65M06)},
  MRNUMBER = {4314112},
MRREVIEWER = {Can Zhang},
      
}

@article {guerrero-2019,
    AUTHOR = {Guerrero, S. and Kassab, K.},
     TITLE = {Carleman estimate and null controllability of a fourth order
              parabolic equation in dimension {$N\geq 2$}},
   JOURNAL = {J. Math. Pures Appl. (9)},
  FJOURNAL = {Journal de Math\'{e}matiques Pures et Appliqu\'{e}es. Neuvi\`eme S\'{e}rie},
    VOLUME = {121},
      YEAR = {2019},
     PAGES = {135--161},
      ISSN = {0021-7824},
   MRCLASS = {93B05 (35K35 93C20)},
  MRNUMBER = {3906168},
MRREVIEWER = {Peng Gao},
       DOI = {10.1016/j.matpur.2018.04.004},
       URL = {},
}

@book {lions,
    AUTHOR = {Lions, J.-L.},
     TITLE = {Contr\^{o}labilit\'{e} exacte, perturbations et stabilisation de
              syst\`emes distribu\'{e}s. {T}ome 1},
    SERIES = {Recherches en Math\'{e}matiques Appliqu\'{e}es [Research in Applied
              Mathematics]},
    VOLUME = {8},
      NOTE = {Contr\^{o}labilit\'{e} exacte. [Exact controllability],
              With appendices by E. Zuazua, C. Bardos, G. Lebeau and J.
              Rauch},
 PUBLISHER = {Masson, Paris},
      YEAR = {1988},
     PAGES = {x+541},
      ISBN = {2-225-81477-5},
   MRCLASS = {49E15 (35B37 49E25 73K10 93B05 93C20)},
  MRNUMBER = {953547},
MRREVIEWER = {Ovidiu C\^{a}rj\u{a}},
}

@article {allonsius,
    AUTHOR = {Allonsius, D. and Boyer, F.},
     TITLE = {Boundary null-controllability of semi-discrete coupled
              parabolic systems in some multi-dimensional geometries},
   JOURNAL = {Math. Control Relat. Fields},
  FJOURNAL = {Mathematical Control and Related Fields},
    VOLUME = {10},
      YEAR = {2020},
    NUMBER = {2},
     PAGES = {217--256},
      ISSN = {2156-8472},
   MRCLASS = {93B05 (35K51 65M06 93-08 93C20)},
  MRNUMBER = {4097728},
MRREVIEWER = {Can Zhang},
       DOI = {10.3934/mcrf.2019037},
       URL = {},
}

@article {allonsius2,
    AUTHOR = {Allonsius, Damien. and Boyer, F. and Morancey, M.},
     TITLE = {Spectral analysis of discrete elliptic operators and
              applications in control theory},
   JOURNAL = {Numer. Math.},
  FJOURNAL = {Numerische Mathematik},
    VOLUME = {140},
      YEAR = {2018},
    NUMBER = {4},
     PAGES = {857--911},
      ISSN = {0029-599X},
   MRCLASS = {34L16 (35K10 65M06 93B05 93B40)},
  MRNUMBER = {3864704},
MRREVIEWER = {Can Zhang},
       DOI = {10.1007/s00211-018-0983-1},
       URL = {},
}

@article {boyer-hernandez,
    AUTHOR = {Boyer, F. and Hern\'{a}ndez-Santamar\'{\i}a, V.},
     TITLE = {Carleman estimates for time-discrete parabolic equations and
              applications to controllability},
   JOURNAL = {ESAIM Control Optim. Calc. Var.},
  FJOURNAL = {ESAIM. Control, Optimisation and Calculus of Variations},
    VOLUME = {26},
      YEAR = {2020},
     PAGES = {Paper No. 12, 43},
      ISSN = {1292-8119},
   MRCLASS = {93B05 (35K20 65M06 93B07 93C55)},
  MRNUMBER = {4064474},
       DOI = {10.1051/cocv/2019072},
       URL = {},
}

@article {cerpa-2010,
    AUTHOR = {Cerpa, E.},
     TITLE = {Null controllability and stabilization of the linear
              {K}uramoto-{S}ivashinsky equation},
   JOURNAL = {Commun. Pure Appl. Anal.},
  FJOURNAL = {Communications on Pure and Applied Analysis},
    VOLUME = {9},
      YEAR = {2010},
    NUMBER = {1},
     PAGES = {91--102},
      ISSN = {1534-0392},
   MRCLASS = {93B05 (35Q53 93C20 93D15)},
  MRNUMBER = {2556747},
MRREVIEWER = {Patrick Martinez},
       DOI = {10.3934/cpaa.2010.9.91},
       URL = {},
}

@article{cahn-1958,
  title={Free energy of a nonuniform system. {I}. {I}nterfacial free energy},
  author={Cahn, John W and Hilliard, John E},
  journal={The Journal of chemical physics},
  volume={28},
  number={2},
  pages={258--267},
  year={1958},
  publisher={American Institute of Physics}
}

@article {LT:2006,
    AUTHOR = {Labb\'{e}, St\'{e}phane and Tr\'{e}lat, Emmanuel},
     TITLE = {Uniform controllability of semidiscrete approximations of
              parabolic control systems},
   JOURNAL = {Systems Control Lett.},
  FJOURNAL = {Systems \& Control Letters},
    VOLUME = {55},
      YEAR = {2006},
    NUMBER = {7},
     PAGES = {597--609},
      ISSN = {0167-6911},
   MRCLASS = {93B05 (35K05 93C25)},
  MRNUMBER = {2225370},
MRREVIEWER = {Pham Tran Nhu},
       DOI = {10.1016/j.sysconle.2006.01.004},
       URL = {},
}

@article {boyer2011,
    AUTHOR = {Boyer, F. and Hubert, F. and Le Rousseau, J.},
     TITLE = {Uniform controllability properties for space/time-discretized
              parabolic equations},
   JOURNAL = {Numer. Math.},
  FJOURNAL = {Numerische Mathematik},
    VOLUME = {118},
      YEAR = {2011},
    NUMBER = {4},
     PAGES = {601--661},
      ISSN = {0029-599X},
   MRCLASS = {93B40 (35K20 65M06 93B05 93B07)},
  MRNUMBER = {2822494},
MRREVIEWER = {Wenming Bian},
       DOI = {10.1007/s00211-011-0368-1},
       URL = {},
}

@article {f-r,
    AUTHOR = {Fattorini, H. O. and Russell, D. L.},
     TITLE = {Exact controllability theorems for linear parabolic equations
              in one space dimension},
   JOURNAL = {Arch. Rational Mech. Anal.},
  FJOURNAL = {Archive for Rational Mechanics and Analysis},
    VOLUME = {43},
      YEAR = {1971},
     PAGES = {272--292},
      ISSN = {0003-9527},
   MRCLASS = {93B05},
  MRNUMBER = {335014},
MRREVIEWER = {F. M. Kirillova},
       DOI = {10.1007/BF00250466},
       URL = {},
}

@article {l-r,
    AUTHOR = {Lebeau, G. and Robbiano, L.},
     TITLE = {Contr\^{o}le exact de l'\'{e}quation de la chaleur},
   JOURNAL = {Comm. Partial Differential Equations},
  FJOURNAL = {Communications in Partial Differential Equations},
    VOLUME = {20},
      YEAR = {1995},
    NUMBER = {1-2},
     PAGES = {335--356},
      ISSN = {0360-5302},
   MRCLASS = {93C20 (35B37 35K05)},
  MRNUMBER = {1312710},
MRREVIEWER = {Fumio Kojima},
       DOI = {10.1080/03605309508821097},
       URL = {},
}

@article {boyer-2014,
    AUTHOR = {Boyer, F. and Le Rousseau, J.},
     TITLE = {Carleman estimates for semi-discrete parabolic operators and
              application to the controllability of semi-linear
              semi-discrete parabolic equations},
   JOURNAL = {Ann. Inst. H. Poincar\'{e} Anal. Non Lin\'{e}aire},
  FJOURNAL = {Annales de l'Institut Henri Poincar\'{e}. Analyse Non Lin\'{e}aire},
    VOLUME = {31},
      YEAR = {2014},
    NUMBER = {5},
     PAGES = {1035--1078},
      ISSN = {0294-1449},
   MRCLASS = {35K58 (35B45 93B05)},
  MRNUMBER = {3258365},
MRREVIEWER = {Joseph L. Shomberg},
       DOI = {10.1016/j.anihpc.2013.07.011},
       URL = {},
}

@article {BHL10a,
    AUTHOR = {Boyer, F. and Hubert, F. and Le Rousseau, J.},
     TITLE = {Discrete {C}arleman estimates for elliptic operators and
              uniform controllability of semi-discretized parabolic
              equations},
   JOURNAL = {J. Math. Pures Appl. (9)},
  FJOURNAL = {Journal de Math\'{e}matiques Pures et Appliqu\'{e}es. Neuvi\`eme S\'{e}rie},
    VOLUME = {93},
      YEAR = {2010},
    NUMBER = {3},
     PAGES = {240--276},
      ISSN = {0021-7824},
   MRCLASS = {93B05 (35K10 65M06 93C20)},
  MRNUMBER = {2601332},
       DOI = {10.1016/j.matpur.2009.11.003},
       URL = {},
}

@article {boyer-sicon,
    AUTHOR = {Boyer, F. and Hubert, F. and Le Rousseau, J.},
     TITLE = {Discrete {C}arleman estimates for elliptic operators in
              arbitrary dimension and applications},
   JOURNAL = {SIAM J. Control Optim.},
  FJOURNAL = {SIAM Journal on Control and Optimization},
    VOLUME = {48},
      YEAR = {2010},
    NUMBER = {8},
     PAGES = {5357--5397},
      ISSN = {0363-0129},
   MRCLASS = {35K20 (35P15 65M06 93B05 93B07)},
  MRNUMBER = {2745778},
       DOI = {10.1137/100784278},
       URL = {},
}

@book {FI,
    AUTHOR = {Fursikov, A. V. and Imanuvilov, O. Yu.},
     TITLE = {Controllability of evolution equations},
    SERIES = {Lecture Notes Series},
    VOLUME = {34},
 PUBLISHER = {Seoul National University, Research Institute of Mathematics,
              Global Analysis Research Center, Seoul},
      YEAR = {1996},
     PAGES = {iv+163},
   MRCLASS = {93-02 (35B37 35Q30 93B05 93C20)},
  MRNUMBER = {1406566},
MRREVIEWER = {Vilmos Komornik},
}

@article {fc-m,
    AUTHOR = {Fern\'{a}ndez-Cara, E. and M\"{u}nch, A.},
     TITLE = {Numerical exact controllability of the 1{D} heat equation:
              duality and {C}arleman weights},
   JOURNAL = {J. Optim. Theory Appl.},
  FJOURNAL = {Journal of Optimization Theory and Applications},
    VOLUME = {163},
      YEAR = {2014},
    NUMBER = {1},
     PAGES = {253--285},
      ISSN = {0022-3239},
   MRCLASS = {93B40 (35K20 65M60 93C20)},
  MRNUMBER = {3260984},
MRREVIEWER = {Sergei A. Ivanov},
       DOI = {10.1007/s10957-013-0517-z},
       URL = {},
}

@article {fc-m2,
    AUTHOR = {Fern\'{a}ndez-Cara, E. and M\"{u}nch, A.},
     TITLE = {Numerical null controllability of semi-linear 1-{D} heat
              equations: fixed point, least squares and {N}ewton methods},
   JOURNAL = {Math. Control Relat. Fields},
  FJOURNAL = {Mathematical Control and Related Fields},
    VOLUME = {2},
      YEAR = {2012},
    NUMBER = {3},
     PAGES = {217--246},
      ISSN = {2156-8472},
   MRCLASS = {93B40 (47H10 49J05 49M15 93B05 93C20)},
  MRNUMBER = {2991568},
MRREVIEWER = {Enrique Zuazua},
       DOI = {10.3934/mcrf.2012.2.217},
       URL = {},
}

@article {c-g-l,
    AUTHOR = {Carthel, C. and Glowinski, R. and Lions, J.-L.},
     TITLE = {On exact and approximate boundary controllabilities for the
              heat equation: a numerical approach},
   JOURNAL = {J. Optim. Theory Appl.},
  FJOURNAL = {Journal of Optimization Theory and Applications},
    VOLUME = {82},
      YEAR = {1994},
    NUMBER = {3},
     PAGES = {429--484},
      ISSN = {0022-3239},
   MRCLASS = {93B05 (65M60 93C20)},
  MRNUMBER = {1290658},
MRREVIEWER = {Sergei A. Ivanov},
       DOI = {10.1007/BF02192213},
       URL = {},
}

@article {m-z,
    AUTHOR = {M\"{u}nch, A. and Zuazua, E.},
     TITLE = {Numerical approximation of null controls for the heat
              equation: ill-posedness and remedies},
   JOURNAL = {Inverse Problems},
  FJOURNAL = {Inverse Problems. An International Journal on the Theory and
              Practice of Inverse Problems, Inverse Methods and Computerized
              Inversion of Data},
    VOLUME = {26},
      YEAR = {2010},
    NUMBER = {8},
     PAGES = {085018, 39},
      ISSN = {0266-5611},
   MRCLASS = {65M32 (35K20 93B40 93C20)},
  MRNUMBER = {2661697},
MRREVIEWER = {Durdimurat K. Durdiev},
       DOI = {10.1088/0266-5611/26/8/085018},
       URL = {},
}

@article {Thuy,
    AUTHOR = {Nguyen, T. N. T.},
     TITLE = {Carleman estimates for semi-discrete parabolic operators with
              a discontinuous diffusion coefficient and applications to
              controllability},
   JOURNAL = {Math. Control Relat. Fields},
  FJOURNAL = {Mathematical Control and Related Fields},
    VOLUME = {4},
      YEAR = {2014},
    NUMBER = {2},
     PAGES = {203--259},
      ISSN = {2156-8472},
   MRCLASS = {93B05 (35K20 65M20 93C20)},
  MRNUMBER = {3167932},
MRREVIEWER = {Chuang Zheng},
       DOI = {10.3934/mcrf.2014.4.203},
       URL = {},
}

@article {Thuy2,
    AUTHOR = {Nguyen, T. N. T.},
     TITLE = {Uniform controllability of semidiscrete approximations for
              parabolic systems in {B}anach spaces},
   JOURNAL = {Discrete Contin. Dyn. Syst. Ser. B},
  FJOURNAL = {Discrete and Continuous Dynamical Systems. Series B. A Journal
              Bridging Mathematics and Sciences},
    VOLUME = {20},
      YEAR = {2015},
    NUMBER = {2},
     PAGES = {613--640},
      ISSN = {1531-3492},
   MRCLASS = {93B05 (65J10 93C20)},
  MRNUMBER = {3331671},
MRREVIEWER = {Zhongcheng Zhou},
       DOI = {10.3934/dcdsb.2015.20.613},
       URL = {},
}

@article {CLTP:2022,
    AUTHOR = {Cerpa, Eduardo and Lecaros, Rodrigo and Nguyen, Thuy N. T. and
              P\'{e}rez, Ariel},
     TITLE = {Carleman estimates and controllability for a semi-discrete
              fourth-order parabolic equation},
   JOURNAL = {J. Math. Pures Appl. (9)},
  FJOURNAL = {Journal de Math\'{e}matiques Pures et Appliqu\'{e}es. Neuvi\`eme S\'{e}rie},
    VOLUME = {164},
      YEAR = {2022},
     PAGES = {93--130},
      ISSN = {0021-7824},
   MRCLASS = {35K52 (65M06 93B05 93B07)},
  MRNUMBER = {4450878},
       DOI = {10.1016/j.matpur.2022.06.003},
       }

@article {cerpa-mercado,
    AUTHOR = {Cerpa, E. and Mercado, A.},
     TITLE = {Local exact controllability to the trajectories of the 1-{D}
              {K}uramoto-{S}ivashinsky equation},
   JOURNAL = {J. Differential Equations},
  FJOURNAL = {Journal of Differential Equations},
    VOLUME = {250},
      YEAR = {2011},
    NUMBER = {4},
     PAGES = {2024--2044},
      ISSN = {0022-0396},
   MRCLASS = {93B05 (35K25 35K55 35Q53 93C20)},
  MRNUMBER = {2763563},
MRREVIEWER = {Patrick Martinez},
       DOI = {10.1016/j.jde.2010.12.015},
       URL = {},
}

@incollection {boyer-canum,
    AUTHOR = {Boyer, F.},
     TITLE = {On the penalised {HUM} approach and its applications to the
              numerical approximation of null-controls for parabolic
              problems},
 BOOKTITLE = {C{ANUM} 2012, 41e {C}ongr\`es {N}ational d'{A}nalyse
              {N}um\'{e}rique},
    SERIES = {ESAIM Proc.},
    VOLUME = {41},
     PAGES = {15--58},
 PUBLISHER = {EDP Sci., Les Ulis},
      YEAR = {2013},
   MRCLASS = {49M25 (65J05)},
  MRNUMBER = {3174955},
MRREVIEWER = {Songting Luo},
       DOI = {10.1051/proc/201341002},
       URL = {},
}

@article {zhou-2012,
    AUTHOR = {Zhou, Z.},
     TITLE = {Observability estimate and null controllability for
              one-dimensional fourth order parabolic equation},
   JOURNAL = {Taiwanese J. Math.},
  FJOURNAL = {Taiwanese Journal of Mathematics},
    VOLUME = {16},
      YEAR = {2012},
    NUMBER = {6},
     PAGES = {1991--2017},
      ISSN = {1027-5487},
   MRCLASS = {93B07 (35K35 93C20)},
  MRNUMBER = {3001831},
       DOI = {10.11650/twjm/1500406835},
       URL = {},
}
\end{document}